\documentclass[final, 1p,times, number]{elsarticle}

\usepackage{algorithm}
\usepackage{algorithmic}
\usepackage{amsmath}
\usepackage{amsthm}

\newproof{pf}{Proof}
\newtheorem{lemma}{Lemma}
\newtheorem{remark}{Remark}
\newtheorem{theorem}{Theorem}
\newtheorem{corollary}{Corollary}
\newtheorem{assumption}{Assumption}
\newtheorem{example}{Example}

\usepackage{color,xcolor}

\usepackage{amssymb}

\begin{document}

\begin{frontmatter}


 \title{Rejuvenating AMLI-Cycle: From Chebyshev Polynomials to Momentum Acceleration
 }

 \author[Chunyan Niu]{Chunyan Niu}
 \ead{chunyanniu@zzu.edu.cn} 
 \author[Yunhui He]{Yunhui He}
 \ead{yhe43@central.uh.edu} 
 \author[Xiaozhe Hu]{Xiaozhe Hu\corref{corresponding}}
 \ead{Xiaozhe.Hu@tufts.edu}
 
 \cortext[corresponding]{Corresponding author}

\address[Chunyan Niu]{School of Mathematics and Statistics, Zhengzhou University, Zhengzhou, 450001, Henan, China}
\address[Yunhui He]{Department of Mathematics, University of Houston,  Houston, 77204-3008, TX, USA}
\address[Xiaozhe Hu]{Department of Mathematics, Tufts University, Medford, 02155,  MA, USA}

\begin{abstract}
In this paper, we investigate the AMLI-cycle method and make two contributions. First, we revisit the AMLI-cycle using the Chebyshev polynomials and establish a theory for its uniform convergence, assuming the two-grid method converges uniformly. This removes the need for estimating extreme eigenvalues at all coarse levels. Only an estimation of the two-grid convergence rate is needed, which could be done on the second coarsest level, simplifying implementation and reducing computational costs for large-scale problems. Second, we introduce a momentum-accelerated AMLI-cycle using polynomials from momentum accelerations. This novel approach ensures a uniform condition number without requiring extreme eigenvalue or two-grid convergence rate estimations, making its implementation as straightforward as standard multigrid methods.  We prove that it is asymptotically as good as the AMLI-cycle using the Chebyshev polynomials when the quadratic momentum-accelerated polynomials is used. Numerical experiments confirm the robustness and efficiency of the momentum-accelerated AMLI-cycle across various problems, demonstrating performance comparable to the Chebyshev-based AMLI-cycle. These findings validate the theoretical advantages and practical efficacy of the momentum-accelerated AMLI-cycle.  \\
\end{abstract}




\end{frontmatter}

\section{Introduction}\label{sec:intro}
Research on multigrid (MG) methods~\cite{Alcouffe81, Dendy82, Dendy83} has been active for decades, driven by their efficiency, scalability, and computational optimality in solving sparse linear systems arising from the discretization of partial differential equations (PDEs). These qualities have led to widespread use of MG methods in practical applications~\cite{Hackbusch85, Bramble93, Briggs00, Trottenberg01, Xu02, Vassilevski08}, particularly the algebraic multigrid (AMG) methods~\cite{Brandt85, Ruge87, Xu17,LinCowenHescottHu2018a}. However, the performance and efficiency of standard V- or W-cycle MG methods can deteriorate as the physical and geometric complexity of the underlying PDEs increases.

For symmetric positive definite (SPD) problems, more advanced multigrid cycles have been developed in the literature. Axelsson and Vassilevski introduced the algebraic multilevel iteration (AMLI) cycle~\cite{Axelsson89, Axelsson90, Vassilevski92}, which employs Chebyshev polynomials to define the coarse-level solver. However, the AMLI-cycle requires accurate estimation of extreme eigenvalues at all coarse levels to compute the polynomial coefficients, which can be challenging in practice. To address this, the K-cycle method~\cite{Axelsson94, Kraus02} was developed as a nonlinear variant of the AMLI-cycle, avoiding the need for eigenvalue estimation through the use of nonlinear preconditioning~\cite{Axelsson91, Golub99, Saad03}.  In the K-cycle MG method, $k$ steps of the nonlinear preconditioned conjugate gradient (NPCG) method are used to define the coarse-level solver, with multigrid on a coarser level acting as a preconditioner. Assuming the convergence factor of the V-cycle multigrid with a bounded-level difference is bounded, the K-cycle was shown to exhibit uniform convergence if $k$ is sufficiently large~\cite{Kraus02}. Additionally, a comparative analysis in~\cite{Hu13} demonstrated that the K-cycle consistently outperforms (or is at least equivalent to) the corresponding $k$-fold V-cycle ($k$V-cycle). However, while the K-cycle eliminates the need for eigenvalue estimation, its nonlinear nature requires the use of the NPCG method, which increases computational and memory costs by losing the three-term recurrence property of the standard conjugate gradient (CG) method.

Our goal is to develop a linear counterpart to the K-cycle, specially, we want to identify suitable polynomials that enable the corresponding linear AMLI-cycle to maintain a uniformly bounded conditioned number without requiring extreme eigenvalue estimations. Therefore, in this work, we first revisit the theoretical analysis of the AMLI-cycle using Chebyshev polynomials and establish a condition for its uniform convergence, assuming uniform convergence of the two-grid method. Our new theoretical results eliminate the need to estimate extreme eigenvalues on all coarse levels, requiring only an estimation of the two-grid convergence rate, which can be determined theoretically or through a two-grid cycle on a coarse level. This simplifies the practical implementation and reduces computational costs, making the AMLI-cycle more efficient for large-scale problems involving many levels.

To further simplify implementation, we explore polynomials derived from momentum acceleration (MA) techniques used in the optimization community, such as Nesterov acceleration (NA) \cite{Nesterov83, Nesterov18} and stationary Anderson Acceleration (sAA(1))
~\cite{sterck2021asymptotic}. We then define the momentum-accelerated AMLI-cycle using these momentum-accelerated polynomials. Theoretically, we establish its uniform condition number estimation under the standard assumption that the two-grid method converges uniformly, without requiring extreme eigenvalue estimation. Additionally, we prove that with quadratic MA polynomials, the momentum-accelerated AMLI-cycle achieves performance asymptotically equivalent to the AMLI-cycle using quadratic Chebyshev polynomials, while higher-order MA polynomials are slightly worse than the corresponding Chebyshev polynomials. Crucially, implementing the momentum-accelerated AMLI-cycle eliminates the need to estimate extreme eigenvalues or the two-grid convergence rate, making it nearly as simple as standard W-cycle or $k$V-cycle MG methods.

Our preliminary numerical tests demonstrate that the AMLI-cycle implementation with Chebyshev polynomials is robust across a wide range of problems, provided a two-grid convergence rate estimation is available at the coarse level. However, consistent with our theoretical analysis, its performance deteriorates in the absence of this estimation. In contrast, the momentum-accelerated AMLI-cycle is universally applicable and exhibits robust and efficient performance in all cases. It achieves results comparable to the Chebyshev-based AMLI-cycle with two-grid convergence rate estimation and the K-cycle in term of iteration counts.  Moreover, in terms of the CPU time, it matches the performance of the Chebyshev-based AMLI-cycle and is approximately twice as fast as the K-cycle across various test problems.  These findings validate our theoretical results and highlight the practical potential of the momentum-accelerated AMLI-cycle.

The rest of this paper is organized as follows. In Section~\ref{sec:AMLI}, we review the AMLI-cycle and its standard theoretical analysis. In Section~\ref{sec:Polynomials from Momentum Acceleration Methods}, we introduce the unified momentum acceleration method for preconditioned linear systems and the corresponding polynomials. Our new theoretical and practical study of the AMLI-cycle using Chebyshev polynomials is presented in Section~\ref{sec:AMLI-Cycle-Chebyshev}. The momentum-accelerated AMLI-cycle methods and their uniformly bounded condition number results are discussed in Section~\ref{sec:AMLI-Cycles-MA}. In Section~\ref{sec:num}, we provide numerical experiments demonstrating the efficiency of the AMLI-cycle with Chebyshev polynomials and momentum-accelerated polynomials. Finally, we present our conclusions in Section~\ref{sec:conc}.

\section{AMLI-cycle}\label{sec:AMLI}
In this section, we review the AMLI-cycle method introduced in~\cite{Axelsson89, Axelsson90, Vassilevski92}, providing a overview of both its algorithmic framework and theoretical development.

We consider solving the following linear system
\begin{equation} \label{eqn:linear-system}
	A \mathbf{x} = \mathbf{b},
\end{equation}
where~$A \in \mathbb{R}^{n \times n}$ is SPD. Assume we have constructed a hierarchical structure of the matrices~$A_\ell \in \mathbb{R}^{n_\ell \times n_\ell}$, $\ell = 1, 2, \cdots,J$, along with~$A_1 = A$, the prolongations~$P_{\ell} \in \mathbb{R}^{n_\ell \times n_{\ell+1}}$, $\ell = 1, 2, \cdots,J-1$, and the restrictions~$R_{\ell} = P_{\ell}^T \in \mathbb{R}^{n_{\ell+1} \times n_{\ell}}$, $\ell = 1, 2, \cdots,J-1$.  Here, we assume that~$A_{\ell+1} = R_{\ell}A_{\ell}P_{\ell}$, $\ell = 1, 2, \cdots,J-1$.  Additionally, let~$M_{\ell}$~denote the smoother on level~$\ell$, such as Jacobi or Gauss-Seidel method.  We define the AMLI-cycle based on a suitably chosen polynomial $p_k(x)$. Algorithm~\ref{alg:AMLIcycle}~summarizes the AMLI-cycle method.
\begin{algorithm}[H]
	\caption{AMLI-cycle MG:~$B_\ell \mathbf{b}$} \label{alg:AMLIcycle}
	\begin{algorithmic}[1]
		\IF{$\ell == J$}
		\STATE $\mathbf{x_\ell} = A_{\ell}^{-1} \mathbf{b}$
		\ELSE
		\STATE $\textbf{Presmoothing}$: $\mathbf{x}_{\ell} \leftarrow  M_{\ell}\mathbf{b}$
		\STATE $\textbf{Coarse-grid\ correction}$: $\mathbf{x}_{\ell} \leftarrow \mathbf{x}_{\ell} + P_{\ell} \widehat{B}_{\ell+1} R_{\ell}(\mathbf{b} - A_\ell\mathbf{x}_{\ell})$,
		where~$\widehat{B}_{\ell+1} = (I - p_k(B_{\ell+1}A_{\ell+1}))A_{\ell+1}^{-1}$
		\STATE $\textbf{Postsmoothing}$: $\mathbf{x}_{\ell} \leftarrow \mathbf{x}_{\ell} + M_{\ell}^{T} (\mathbf{b} - A_\ell\mathbf{x}_{\ell})$
		\ENDIF
		\STATE $B_{\ell} \mathbf{b} \gets \mathbf{x}_{\ell}$
	\end{algorithmic}
\end{algorithm}

From Algorithm~\ref{alg:AMLIcycle}, the AMLI-cycle defines the operators $B_{\ell}$, $\ell = 1, 2, \cdots,$ $ J-1$,
\begin{align}\label{def:B}
	B_{\ell}=\overline{M}_{\ell}+(I-M^{T}_{\ell}A_{\ell})P_{\ell}\widehat{B}_{\ell+1}P^{T}_{\ell}(I-A_{\ell}M_{\ell}),
\end{align}
where $\overline{M}_{\ell}= M_{\ell}(M^{-1}_{\ell}+M^{-T}_{\ell}-A)M^{T}_{\ell}$ is the symmetrization of $M_{\ell}$ and $\widehat{B}_J = A_J^{-1}$. \eqref{def:B} is a useful identity for the theoretical analysis of the AMLI-cycle. 

The basic requirement of the polynomial is that $p_k(0) = 1$, Thus, we can also define another polynomial $ \displaystyle q_{k-1}(x) := \frac{1-p_k(x)}{x}. $ This implies that
\begin{equation}\label{def:B_hat}
	\widehat{B}_{\ell} = (I - p_k(B_{\ell}A_{\ell}))A_{\ell}^{-1} = B_{\ell} q_{k-1}(A_{\ell} B_{\ell}).
\end{equation}
Based on~\eqref{def:B_hat}, we have the following lemma which shows that the positivity of $q_{k-1}(x)$ implies that $\widehat{B}_{\ell}$ is positive definite.
\begin{lemma}\label{lem:B_hat_SPD}
	Assume that $A_{\ell}$ and $B_{\ell}$ are SPD.  If $q_{k-1}(x) > 0$, for~$0<\lambda_{\min}(B_{\ell}A_{\ell}) \leq x \leq \lambda_{\max}(B_{\ell}A_{\ell})$, or equivalently, $p_{k}(x) < 1$, for~$0<\lambda_{\min}(B_{\ell}A_{\ell}) $ $\leq x \leq \lambda_{\max}(B_{\ell}A_{\ell})$, then $\widehat{B}_{\ell}$ is SPD.
\end{lemma}
\begin{pf}
	By~\eqref{def:B_hat}, $\widehat{B}_{\ell}$ is symmetric.  Furthermore, for any $\mathbf{v}$, we have, 
	\begin{align*}
		(\widehat{B}_{\ell} \mathbf{v}, \mathbf{v}) = (B_{\ell} q_{k-1}(A_{\ell}B_{\ell})\mathbf{v}, \mathbf{v}) = (q_{k-1} (B_{\ell}^{1/2} A_{\ell} B_{\ell}^{1/2})B_{\ell}^{1/2}\mathbf{v}, B_{\ell}^{1/2}\mathbf{v}).
	\end{align*}
	Note that $\lambda(B_{\ell}^{1/2} A_{\ell} B_{\ell}^{1/2}) = \lambda(B_{\ell}A_{\ell})$, we have
	\begin{align*}
		(\widehat{B}_{\ell} \mathbf{v}, \mathbf{v})  \geq \min_{x \in [\lambda_{\min}(B_{\ell}A_{\ell}), \lambda_{\max}(B_{\ell}A_{\ell})]} q_{k-1}(x) (B_{\ell} \mathbf{v}, \mathbf{v}) > 0,
	\end{align*}
	if $q_{k-1}(x) > 0$ for $0<\lambda_{\min}(B_{\ell}A_{\ell}) \leq x \leq \lambda_{\max}(B_{\ell}A_{\ell})$.  This implies that $\widehat{B}_{\ell}$ is SPD and completes the proof.
\end{pf}

Next lemma shows that if $\widehat{B}_{\ell+1}$ is SPD, so is $B_{\ell}$, provided that $M_{\ell}$ is nonexpansive. 
\begin{lemma}\label{lem:B_SPD}
	If $\| I - M_{\ell} A_{\ell} \|_{A_{\ell}} \leq  1$ and $\widehat{B}_{\ell+1}$ is SPD, then $B_{\ell}$ is SPD.
\end{lemma}
\begin{pf}
	From~\eqref{def:B}, it is easy to see that $B_{\ell}$ is symmetric.  For any $\mathbf{v} \neq 0$, we have
	\begin{align*}
		((I-B_{\ell}A_{\ell})\mathbf{v}, \mathbf{v})_{A_{\ell}} &= ( (I - P_{\ell} \widehat{B}_{\ell+1} R_{\ell} A_{\ell})\widetilde{\mathbf{v}}, \widetilde{\mathbf{v}} )_{A_{\ell}} = \| \widetilde{\mathbf{v}} \|_{A_{\ell}}^2 - (R_{\ell}A_{\ell} \widetilde{\mathbf{v}},  R_{\ell}A_{\ell} \widetilde{\mathbf{v}})_{\widehat{B}_{\ell+1}},
	\end{align*}
	where $\widetilde{\mathbf{v}}:= (I - M_{\ell} A_{\ell}) \mathbf{v}$.  Since $\widehat{B}_{\ell+1}$ is SPD, if $\widetilde{\mathbf{v}} \neq 0$, we can get
	$\| \mathbf{v} \|^2_{A_{\ell}} - (B_{\ell} A_{\ell} \mathbf{v}, \mathbf{v})_{A_{\ell}} < \| \widetilde{\mathbf{v}} \|^2_{A_{\ell}} \leq \| \mathbf{v} \|^2_{A_{\ell}}$.
	If  $\widetilde{\mathbf{v}} = 0$, we have
	$	\| \mathbf{v} \|^2_{A_{\ell}} - (B_{\ell} A_{\ell} \mathbf{v}, \mathbf{v})_{A_{\ell}} = 0 < \| \mathbf{v} \|^2_{A_{\ell}}$.
	Thus, $(B_{\ell} A_{\ell} \mathbf{v}, \mathbf{v})_{A_{\ell}} >0 $ which means that $B_{\ell}$ is SPD.  This completes the proof.
\end{pf}

Based on Lemma~\ref{lem:B_hat_SPD} and Lemma~\ref{lem:B_SPD}, we immediately see that the operator $B_{\ell}$, $\ell=1, 2, \cdots, J$, defined by the AMLI-cycle, i.e., Algorithm~\ref{alg:AMLIcycle} is SPD. 
\begin{theorem}\label{thm:AMLI-SPD}
	Let the operator $B_{\ell}$, $\ell=1,2, \cdots, J-1$, be defined in~\eqref{def:B} and the operator $\widehat{B}_{\ell}$, $\ell = 2, 3, \cdots, J-1$ be defined in~\eqref{def:B_hat} with a polynomial $p_k(x)$ such that $p_k(0)=1$ and $\widehat{B}_{J} = A_{J}^{-1}$.  Assume that $A_{\ell}$, $\ell = 1, 2, \cdots, J$, are SPD and $p_k(x) <1$ for $x \in [\lambda_{\min}(B_{\ell}A_{\ell}), \lambda_{\max}(B_{\ell}A_{\ell})]$.  We have that the operators $B_{\ell}$, $\ell = 1, ,2, \cdots, J-1$, and $\widehat{B}_{\ell}$, $\ell=2,3,\cdots, J$, are all SPD. 
\end{theorem}
\begin{pf}
	The conclusion follows from the fact that $\widehat{B}_{J} = A_{J}^{-1}$, which is SPD, combined with the recursive application of Lemma~\ref{lem:B_hat_SPD} and Lemma~\ref{lem:B_SPD}.
\end{pf}

The convergence analysis of the AMLI-cycle has been extensively studied in the literature~\cite{Vassilevski08,Axelsson89, Axelsson90, Kraus02}. Here, we revisit a version originally presented in~\cite{Kraus12} and briefly include the proof for completeness.
\begin{theorem} \label{thm:AMLI-convergence}
	Let~$B_{\ell}$~and~$\widehat{B}_{\ell}$~be defined as in~\eqref{def:B} and~\eqref{def:B_hat}, respectively. Under the assumptions of~Theorem~\ref{thm:AMLI-SPD}, the AMLI operator $B_{\ell}$ satisfies the following condition number estimate
	\begin{equation} \label{ine:B_cond_est}
		\kappa(B_{\ell}A_{\ell}) \leq \kappa({B^{TG}_{\ell}A_{\ell}})\frac{\max\{1,\gamma^{\ell+1}_{1}\}}{\min\{1,\gamma^{\ell+1}_{0}\}},
	\end{equation}
	where, for $\mu_{\ell+1} = \lambda_{\min}(B_{\ell+1}A_{\ell+1})$ and~$L_{\ell+1} = \lambda_{\max}(B_{\ell+1}A_{\ell+1})$, 
	\begin{align}
		\gamma^{\ell+1}_{1} = \max_{x \in [\mu_{\ell+1}, L_{\ell+1}]}[xq_{k-1}(x)] = 1-  \min_{x \in [\mu_{\ell+1}, L_{\ell+1}]}[p_{k}(x)] >0, \label{def:r_1}\\
		\gamma^{\ell+1}_{0} = \min_{x \in [\mu_{\ell+1}, L_{\ell+1}]}[xq_{k-1}(x)] = 1-  \max_{x \in [\mu_{\ell+1}, L_{\ell+1}]}[p_{k}(x)] >0. \label{def:r_0}
	\end{align}
	Furthermore, if $0<\mu \leq \mu_{\ell}\leq L_{\ell} \leq L <\infty$ for all levels $\ell$, $\kappa(B^{TG}_{\ell}A_{\ell}) \leq \kappa_{TG}$, and $p_k(x) < 1$, for $x \in [\mu, L]$, then the condition number $\kappa(B_{\ell}A_{\ell})$ is uniformly bounded,
	\begin{align} \label{ine:B_uniform_bound}
		\kappa(B_{\ell}A_{\ell}) \leq \kappa_{TG} \frac{\max\{1, \gamma_1\}}{\min\{1, \gamma_0\}},
	\end{align} 
	where 
	$
	\gamma_1 := 1 - \min_{x \in [\mu, L]} [p_k(x)]$ and
	$\gamma_0 := 1 - \max_{x \in [\mu, L]} [p_k(x)]$.  
\end{theorem}

\begin{pf}
	From~\eqref{def:B}, we have
	\begin{align*}
		(B_{\ell}A_{\ell} \mathbf{v}, \mathbf{v})_{A_{\ell}} 
		& = (\overline{M}_{\ell} A_{\ell} \mathbf{v}, \mathbf{v})_{A_{\ell}} + (A_{\ell+1}^{-1} \mathbf{w}, \mathbf{w}) \frac{(\widehat{B}_{\ell+1}\mathbf{w}, \mathbf{w})}{(A_{\ell+1}^{-1} \mathbf{w}, \mathbf{w})},
	\end{align*}
	where $\mathbf{w} = R_{\ell}A_{\ell}(I-M_{\ell}A_{\ell})\mathbf{v}$. Note that, 
	\begin{align*}
		(\widehat{B}_{\ell+1} \mathbf{w}, \mathbf{w}) &= (B_{\ell+1} q_{k-1}(A_{\ell+1}B_{\ell+1})\mathbf{w}, \mathbf{w}) \\
		&= (A_{\ell+1}^{1/2} B_{\ell+1}A_{\ell+1}^{1/2}q_{k-1} (A_{\ell+1}^{1/2}B_{\ell+1} A_{\ell+1}^{1/2}) A_{\ell+1}^{-1/2}\mathbf{w}, A_{\ell+1}^{-1/2}\mathbf{w}).
	\end{align*}
	Using the definition of $\gamma_0^{\ell+1}$ \eqref{def:r_0} and $\gamma_1^{\ell+1}$ \eqref{def:r_1}, we get 
	\begin{align*}
		\gamma_{0}^{\ell+1} (A_{\ell+1}^{-1} \mathbf{w}, \mathbf{w}) \leq (\widehat{B}_{\ell+1} \mathbf{w}, \mathbf{w}) \leq \gamma_1^{\ell+1} (A_{\ell+1}^{-1} \mathbf{w}, \mathbf{w}).
	\end{align*}
	Thus, 
	\begin{align*}
		(\overline{M}_{\ell} A_{\ell} \mathbf{v}, \mathbf{v})_{A_{\ell}} \hskip -1pt + \hskip -1pt  \gamma_0^{\ell+1} (A_{\ell+1}^{-1}\mathbf{w}, \mathbf{w}) &\leq (B_{\ell}A_{\ell} \mathbf{v}, \mathbf{v})_{A_{\ell}} \\
		& \leq  (\overline{M}_{\ell} A_{\ell} \mathbf{v}, \mathbf{v})_{A_{\ell}} \hskip -1pt + \hskip -1pt  \gamma_1^{\ell+1} (A_{\ell+1}^{-1}\mathbf{w}, \mathbf{w}).
	\end{align*} 
	Since the two-grid operator on level $\ell$ is defined as $B^{TG}_{\ell} := \overline{M}_{\ell}+(I-M^{T}_{\ell}A_{\ell})P_{\ell}A_{\ell+1}^{-1}P^{T}_{\ell}(I-A_{\ell}M_{\ell})$, we arrive at
	\begin{align*}
		\min\{1,\gamma_0^{\ell+1}\} (B^{TG}_{\ell}A_{\ell} \mathbf{v}, \mathbf{v})_{A_{\ell}}  \leq (B_{\ell}A_{\ell} \mathbf{v}, \mathbf{v})_{A_{\ell}}  \leq  \max \{ 1, \gamma_1^{\ell+1} \} (B^{TG}_{\ell}A_{\ell} \mathbf{v}, \mathbf{v})_{A_{\ell}}.
	\end{align*}
	Using the property that $\lambda_{\max}(B^{TG}_{\ell}A_\ell) \leq 1$ and $\lambda_{\min}(B^{TG}_{\ell}A_{\ell}) \geq (1 - \delta^{TG}_{\ell})$, where $\delta^{TG}_\ell$ is the convergence rate of the two-grid method on level $\ell$, we obtain
	\begin{align*}
		\min\{1,\gamma_0^{\ell+1}\} (1 - \delta^{TG}_{\ell}) ( \mathbf{v}, \mathbf{v})_{A_{\ell}}  \leq (B_{\ell}A_{\ell} \mathbf{v}, \mathbf{v})_{A_{\ell}}  \leq  \max \{ 1, \gamma_1^{\ell+1} \} ( \mathbf{v}, \mathbf{v})_{A_{\ell}}.
	\end{align*}
	Thus, \eqref{ine:B_cond_est} follows with
	\begin{align} \label{def:mu_L_recursive}
		\mu_{\ell} =  \min\{1,\gamma_0^{\ell+1}\} (1 - \delta^{TG}_{\ell}) \ \text{and} \  L_{\ell} =  \max \{ 1, \gamma_1^{\ell+1} \}.
	\end{align}

	If $0<\mu \leq \mu_{\ell}\leq L_{\ell} \leq L <\infty$ for all levels $\ell$ and $p_k(x) < 1$, for $x \in [\mu, L]$, we have, for all $\ell$,
	$\gamma_1^{\ell} \leq \gamma_1 := 1 - \min_{x \in [\mu, L]} [p_k(x)]$ and 
	$\gamma_0^{\ell} \geq \gamma_0 := 1 - \max_{x \in [\mu, L]} [p_k(x)]$.  
	Then, if $\kappa(B^{TG}_{\ell}A_{\ell}) \leq \kappa_{TG}$, we get the uniform bound~\eqref{ine:B_uniform_bound} directly from~\eqref{ine:B_cond_est}.
\end{pf}

A simple choice of the polynomial is $p_k(x) = \left(1-\frac{x}{L_{\ell+1}}\right)^k$.  When $k=1$ and $2$, the resulting AMLI-cycle is nothing but the standard V-cycle and W-cycle methods (with $L_{\ell} = 1$, as we will see later).  For general $k$, the corresponding AMLI-cycle is referred to as  the $k$V-cycle.  Note that $0 \leq p_k(x) \leq   \left(1 - \frac{\mu_{\ell+1}}{L_{\ell+1}} \right)^k  < 1$ for $0< \mu_{\ell+1} \leq x \leq L_{\ell+1}$. From~\eqref{def:mu_L_recursive}, it follows that $L_{\ell} \leq 1$ and $\mu_{\ell} = (1 - (1- \frac{\mu_{\ell+1}}{L_{\ell+1}})^k)(1 - \delta^{TG}_{\ell}) \geq 0$.  This allows us to apply the same argument recursively.  Together with the fact that $\mu_{J} = L_{J} = 1$,  we can set the uniform upper bound of $L_{\ell}$ be $L=1$ and,  under the assumption that the two-grid method is uniformly convergent on all levels, i.e., $\delta^{TG}_{\ell} \leq \delta_{TG} < 1$, the uniform lower bound of $\mu_{\ell}$ can be determined by $\mu \leq (1 - (1-\mu)^k) (1 - \delta_{TG})$.   This leads to the following corollary. 

\begin{corollary}
	Assuming the two-grid method converges uniformly with convergence rate $\delta_{TG}$. The AMLI operator $B_{\ell}$, defined by Algorithm~\ref{alg:AMLIcycle} using the simple polynomial $p_k(x) = (1-x)^k$, $k \geq 2$, satisfies the condition number estimate $\kappa(B_{\ell}A_{\ell}) \leq 1/\mu$
	provided there exists $0 < \mu < 1$ satisfying $\mu \leq (1 - (1-\mu)^k) (1 - \delta_{TG})$. Moreover, a sufficient condition of the existence of $\mu$ is 
	\begin{equation} \label{ine:TG-rate-bound-simple}
		\delta_{TG} < 1 - \frac{1}{k}.
	\end{equation}
\end{corollary}
\begin{pf}
	The result directly follows from the observation that we can simple use $L_{\ell} = 1$ on all levels to define the simply polynomial and the following elementary inequality for $\mu > 0$
	\begin{align*}
		\delta_{TG} \leq 1 - \frac{\mu}{ 1- (1-\mu)^k} < 1 - \frac{1}{k}.
	\end{align*}
\end{pf}

\begin{remark}
	When $k=2$, the AMLI-cycle corresponds to the W-cycle and requiring $\delta_{TG} < 0.5$.  This aligns with the well-known result that, when the two-grid method is uniformly convergent with convergence rate less than $0.5$, the multilevel W-cycle also converges uniformly.  See,  for example, \cite{Vassilevski08}.  Notably, with this simple choice of polynomial, there is no need to estimate the extreme eigenvalues, as we simply set $L_{\ell} = 1$ on all levels.
\end{remark} 

Due to its min-max property, the Chebyshev polynomials have been natural choices since the AMLI-cycle was proposed.  As outlined in~\cite{Vassilevski08,Axelsson89,Axelsson90}, the following scaled and shifted Chebyshev polynomials are commonly employed in the AMLI-cycle, 
\begin{align}\label{def:chebyshev}
	p_k(x) = \frac{1 + T_k(\frac{L_\ell+\mu_\ell - 2x}{L_\ell-\mu_\ell})}{1 + T_k(\frac{L_\ell+\mu_\ell}{L_\ell-\mu_\ell})},
\end{align}	
where $T_k(x)$ is the usual Chebyshev polynomials of the first kind. It is easy to verify that $p_k(0)=1$ and $p_k(x) \in [0,1]$ for $x \in [\mu_{\ell},L_{\ell}]$.   As discussed in the literature, e.g., \cite{Vassilevski08,Axelsson89, Axelsson90}, the AMLI-cycle using Chebyshev polynomials converges faster than its $k$V-cycle counterpart. The details of this comparison will be explored later. However, implementing the AMLI-cycle requires estimates of the smallest eigenvalue $\mu_\ell = \lambda_{\min}(B_{\ell}A_{\ell})$ and the largest eigenvalue $L_\ell = \lambda_{\max}(B_{\ell}A_{\ell})$. In general, the performance of the AMLI-cycle, Algorithm~\ref{alg:AMLIcycle}, depends on accurately estimating these extreme eigenvalues. While a good estimate for the largest eigenvalue $L_\ell$ is usually feasible in the SPD case, obtaining a reliable estimate for the smallest eigenvalue $\mu_\ell$ is challenging. 
In~Section~\ref{sec:AMLI-Cycle-Chebyshev}, we present a new theoretical framework for the AMLI-cycle using Chebyshev polynomials, which simplifies its practical implementation.

\begin{remark}
	Other polynomials have also been developed for use in AMLI-cycles. For instance, the polynomial that provides the best approximation of $1/x$ with respect to the uniform norm \cite{Kraus12}, and the polynomial designed for solving graph Laplacians based on matching \cite{brannick2013algebraic}. Like the Chebyshev polynomials, the AMLI-cycle using these alternatives also achieves faster convergence compared to their $k$V-cycle counterparts. However, similar to the Chebyshev polynomials, their implementation still requires estimations of the extreme eigenvalues $\mu_\ell$ and $L_\ell$. This work focuses Chebyshev polynomials, while a comparison with other polynomials is the subject of ongoing research and will be addressed in future publications.
\end{remark}

\section{Polynomials from Momentum Acceleration Methods}\label{sec:Polynomials from Momentum Acceleration Methods}
In this section, we introduce the polynomials that are derived from the momentum acceleration methods originally developed for the optimization problems. 

\subsection{Momentum Acceleration Methods}\label{subsec:Momentum Acceleration Methods}
We begin by introducing momentum acceleration techniques, focusing on first-order methods for addressing the following unconstrained optimization problem:
\begin{equation} \label{eqn:opt problem}
	\min_{\mathbf{x} \in \mathbb{R}^{n}}\mathbf{f}(\mathbf{x}),
\end{equation}
where~$\mathbf{f}:\mathbb{R}^{n}\rightarrow \mathbb{R}$~is a continuously differentiable strongly convex function satisfying
\begin{equation}\label{property of f}
	\frac{\mu}{2}\|\mathbf{x} - \mathbf{y}\|^{2} \leq  \mathbf{f}(\mathbf{x}) - \mathbf{f}(\mathbf{y}) - ( \nabla \mathbf{f}(\mathbf{y}), \mathbf{x} - \mathbf{y}) \leq \frac{L}{2}\|\mathbf{x} - \mathbf{y}\|^{2},\ \forall\  \mathbf{x},\mathbf{y}\in \mathbb{R}^{n},
\end{equation}
with~$L >0$~and~$\mu > 0$~being the Lipschitz and convexity constants, respectively.  Here~$(\cdot,\cdot)$ represents a generic inner product of~$\mathbb{R}^n$~and~$\| \cdot \|$~denotes its induced norm.

The optimization problem~\eqref{eqn:opt problem}~is typically solved by the gradient descent~(GD)~method, which, under suitable assumptions, converges linearly with a rate~$\frac{L-\mu}{L+\mu}$~\cite{Nesterov18}. Various algorithms have been developed to accelerate the convergence, and the main idea is to incorporate some momentum. Next, we present the unified momentum-accelerated GD method in Algorithm~\ref{alg:Unified-MA}.  Different choices of the parameters $\alpha$ and $\beta$ lead to different momentum and result in various methods. We present two choices from the literature in~Table~\ref{tab:choices-MA}.
\begin{algorithm}[h]
	\caption{Unified Momentum Acceleration Method} \label{alg:Unified-MA}
	\begin{algorithmic}[1]
		\STATE $\mathbf{x}^{0}, \mathbf{y}^{0}$~are given as initial iterates. $\alpha$ and $\beta$~are given parameters.
		\FOR{$k =0,1,2,...$}
		\STATE $\mathbf{y}^{k+1} \gets \mathbf{x}^{k} - \alpha \nabla \mathbf{f}(\mathbf{x}^{k})$
		\STATE $\mathbf{x}^{k+1} \gets \mathbf{y}^{k+1} + \beta(\mathbf{y}^{k+1} - \mathbf{y}^k )$
		\ENDFOR
	\end{algorithmic}
\end{algorithm}

\begin{table}[h!]
	\begin{center}
		\caption{Different momentum acceleration methods} \label{tab:choices-MA}
		\begin{tabular}{| c || c | c |}
			\hline \hline 
			&  $\alpha$  & $\beta$ \\
			\hline \hline 
			Nesterov acceleration (NA)~\cite{Nesterov83} & $\displaystyle \frac{1}{L}$ & $\displaystyle\frac{\sqrt{L} - \sqrt{\mu}}{\sqrt{L} + \sqrt{\mu}}$ \\ \hline 
			Stationary Anderson Acceleration (sAA(1))~\cite{sterck2021asymptotic} &  $\displaystyle \frac{4}{\mu + 3 L}$ & $\displaystyle\frac{1 - \sqrt{\alpha\mu}}{1 + \sqrt{\alpha\mu}}$ \\
			\hline  \hline 
		\end{tabular}
	\end{center}	
\end{table}

\begin{remark}
	Replacing $\mathbf{y}^k$ in Algorithm~\ref{alg:Unified-MA} by $\mathbf{x}^{k-1}$ and setting $\alpha = \frac{2}{L+\mu}$ and $\beta = \left( \frac{\sqrt{L} - \sqrt{\mu}}{\sqrt{L}  + \sqrt{\mu}}  \right)^2$, we recover the well-known heavy ball method~\cite{Polyak64}, another prominent momentum acceleration technique.  In this work, however, we focus on the NA and sAA(1) methods. 
\end{remark}

\subsection{Polynomials for Momentum Acceleration}\label{subsec:Polynomials for Momentum Acceleration}
To derive the polynomials that correspond to the momentum acceleration methods~Algorithm~\ref{alg:Unified-MA}, we reformulate the task of solving~\eqref{eqn:linear-system}~with an SPD preconditioner~$B$~as solving the following quadratic optimization problem,
\begin{equation}\label{eqn:opt problem with B}
	\min_{\mathbf{x} \in \mathbb{R}^{n}}\frac{1}{2}(BA \mathbf{x},\mathbf{x})_{B^{-1}} - (B\mathbf{b},\mathbf{x})_{B^{-1}},
\end{equation}
where~$(\mathbf{x},\mathbf{x})_{B^{-1}} := (B^{-1} \mathbf{x}, \mathbf{x})$.  
Applying the MA method~Algorithm~\ref{alg:Unified-MA} to~\eqref{eqn:opt problem with B} and eliminating~$\mathbf{y}^k$~leads to~Algorithm~\ref{alg:MAwithB} as follows
\begin{algorithm}[H]
	\caption{MA method for preconditioned linear systems:~$\widehat{B}^{N}$} \label{alg:MAwithB}
	\begin{algorithmic}[1]
		\STATE $\mathbf{x}^{0}, \mathbf{x}^{1}$~are given as initial iterates and $\alpha$ and $\beta$ are a given parameters.
		\FOR{$i = 2,...,k$}
		\STATE $\mathbf{x}^{i} \gets (1+\beta)[\mathbf{x}^{i-1} + \alpha B(\mathbf{b} - A\mathbf{x}^{i-1})]- \beta[\mathbf{x}^{i-2} + \alpha B(\mathbf{b}- A\mathbf{x}^{i-2})]$
		\ENDFOR
	\end{algorithmic}
\end{algorithm}

Algorithm~\ref{alg:MAwithB}~shows that the MA method for solving preconditioned linear systems is essentially a linear iterative method utilizing the previous two steps.  It computes each update as a weighted average of the last two updates.  To derive the polynomial corresponding to the MA method, we define the error at the $k$-th step of the MA method as~$\mathbf{e}^k := \mathbf{x}^* - \mathbf{x}^k$. From~Algorithm~\ref{alg:MAwithB}, it satisfies the following three-term recurrence relationship:
\begin{equation*} 
	\mathbf{e}^{k+1} = (1+\beta)(I - \alpha BA)\mathbf{e}^{k} - \beta(I - \alpha BA)\mathbf{e}^{k-1}.
\end{equation*}
Choosing $\mathbf{e}^1 = p_1(BA) \mathbf{e}^0$, for example, $\mathbf{e}^1 = (I - \alpha BA) \mathbf{e}^0$ which corresponds to $\mathbf{x}^1 = \mathbf{x}^0 + \alpha B(\mathbf{b} - A \mathbf{x}^{0})$, this implies that~$\mathbf{e}^k = p_k(BA)\mathbf{e}^0$~where~$p_k(x)$~is a polynomial of degree at most~$k$~and satisfies~$p_k(0)=1$,
\begin{equation}\label{eqn:recurrence-relation-of-pk-MA}
	p_{k+1}(x)  = (1+\beta)(1 - \alpha x )p_{k}(x) - \beta(1 - \alpha x)p_{k-1}(x).
\end{equation}
If we choose $\alpha$ and $\beta$ as in~Table~\ref{tab:choices-MA} and substitute them back into~\eqref{eqn:recurrence-relation-of-pk-MA}, we can obtain the MA polynomials for the NA and sAA(1) methods, respectively.  In Section~\ref{sec:AMLI-Cycles-MA}, we use the MA polynomials to define the AMLI-cycle and discuss the choices of $\alpha$ and $\beta$ accordingly.

\section{AMLI-Cycle using the Chebyshev Polynomials} \label{sec:AMLI-Cycle-Chebyshev}
In this section, we present new theoretical advancements of the AMLI-cycle using the Chebyshev polynomial~\eqref{def:chebyshev}, enabling more efficient practical implementation.  Previous studies~\cite{Vassilevski08,Axelsson89,Axelsson90,Vassilevski92,Kraus02} often required accurate estimations of the extreme eigenvalues at each level, which can be computationally costly. In this work, we show that we only need an estimation (upper bound) of the two-grid convergence rate to define the Chebyshev polynomials and the corresponding AMLI-cycle.  Such an estimation could be derived theoretically, as shown in ~\cite{brannick2013algebraic,napov2012algebraic}, or efficiently computed on a coarse level, e.g., the second coarsest level.  Furthermore, we prove the resulting AMLI-cycle achieves uniform convergence with a larger two-grid convergence rate compared with the $k$V-cycles for all $k \geq 2$. This result, to the best of our knowledge, extends the existing result for $k=2$ case~\cite{Vassilevski08} to general $k$.  

Consider the scaled and shifted Chebyshev polynomial~\eqref{def:chebyshev}, defined using the estimation of the extreme eigenvalues $\mu_{\ell+1}$ and $L_{\ell+1}$.  It follows that $ 0 \leq p_k(x) \leq p_k(\mu_{\ell+1}) < 1$ for $0 < \mu_{\ell+1} \leq x \leq L_{\ell+1}$.  Therefore, from~\eqref{def:mu_L_recursive}, we have $L_{\ell} \leq  1$ and $\mu_{\ell} = \left[ 1 -  p_{k}(\mu_{\ell+1}) \right](1 - \delta^{TG}_{\ell}) > 0$, which allows us to do the recursion on all levels.  Since $L_J = 1$, we can set the uniform upper bound of $L_{\ell}$ as $L=1$.  Additionally, this motivates us to derive the uniform lower bound of $\mu_\ell$ by solving $ \mu \leq   \left[ 1 - p_k(\mu) \right](1 - \delta_{TG})$, where $\delta_{TG}$ represents the uniform upper bound of the two-grid method across on all levels.   This leads the following theorem.

\begin{theorem} \label{thm:uniform-conv-AMLI-chebyshev}
	Assuming the two-grid method converges uniformly with convergence rate $\delta_{TG}$. Consider the following scaled and shifted Chebyshev polynomial for $k \geq 2$,
	\begin{align}\label{def:chebyshev-mu}
		p_k(x) = \frac{1 + T_k(\frac{1+\mu - 2x}{1-\mu})}{1 + T_k(\frac{1+\mu}{1-\mu})},
	\end{align}
	where $0< \mu < 1$ satisfies 
	\begin{equation}\label{ine:chebyshev-condition_mu}
		\mu \leq   \left[ 1 - p_k(\mu) \right](1 - \delta_{TG}). 
	\end{equation}
	Then the AMLI operator $B_{\ell}$, defined by Algorithm~\ref{alg:AMLIcycle} using the scaled and shifted Chebyshev polynomial~\eqref{def:chebyshev-mu}, satisfies the condition number estimate $\kappa(B_{\ell} A_{\ell}) \leq 1/\mu$.
	Moreover, a sufficient condition of \eqref{ine:chebyshev-condition_mu} is
	\begin{equation}\label{ine:TG-rate-bound-chebyshev}
		\delta_{TG} < 1 - \frac{1}{k^2}.
	\end{equation}
\end{theorem}

\begin{pf}
	On the coarsest level, we have $B_J = A_J^{-1}$. Thus $0 < \mu \leq \mu_J = L_J = 1$.  On the second coarsest level,  we have $B_{J-1} = B_{J-1}^{TG}$. Therefore, $0 < \mu \leq \mu_{J-1} = 1 - \delta_{TG} < L_{J-1} = 1$. On level $\ell = J-2$, by~\eqref{def:mu_L_recursive},  we have 
	\begin{align*}
		L_{J-2} &= \max\{1, \gamma_1^{J-1} \} = \max \{1, 1 - \min_{x \in [\mu_{J-1}, L_{J-1}]} [p_k(x)]  \} \\
		&\leq \max \{1, 1 - \min_{x \in [\mu, 1]} [p_k(x)]  \} = 1, \\
		\mu_{J-2} &= (1-\delta_{TG}) \min\{ 1, \gamma_0^{J-1} \} = (1-\delta_{TG}) \min\{ 1,  1 - \max_{x \in [\mu_{J-1}, L_{J-1}]} [p_k(x)]   \} \\
		& \geq (1-\delta_{TG}) \min\{ 1,  1 - \max_{x \in [\mu, 1]} [p_k(x)]   \} = (1-\delta_{TG}) [1 - p_k(\mu)].
	\end{align*}
	This implies $0 < \mu \leq \mu_{J-2} < L_{J-2} \leq 1$ by \eqref{ine:chebyshev-condition_mu}, allowing us to apply the same argument recursively.  Therefore, based on Theorem~\ref{thm:AMLI-convergence}, we obtain the condition number estimate.
	
	To show the sufficient condition \eqref{ine:TG-rate-bound-chebyshev}, from \eqref{ine:chebyshev-condition_mu}, we have  
	\begin{equation*}
		\delta_{TG} \leq  1 - \frac{\mu}{1- p_k(\mu)} < 1 - \lim_{ \mu \mapsto 0^+} \frac{\mu}{ 1 - p_{k}(\mu)}.
	\end{equation*}
	By the L'H\^{o}pital's rule, we have
	$	\lim_{ \mu \mapsto 0^+} \frac{\mu}{ 1 - p_{k}(\mu)} = \lim_{ \mu \mapsto 0^+} \frac{1}{- p_k'(\mu)}$.
	Note that, using the property $T_k'(x) = k U_{k-1}(x)$ where $U_{k}(x)$ is the Chebyshev polynomial of the second kind, 
	\begin{align*}
		p'_k(\mu) &= \left( \frac{2}{1 + T_{k}(\frac{1+\mu}{1-\mu})}  \right)'\\
        &= (-2) \left(1 +  T_{k}\left(\frac{1+\mu}{1-\mu}\right) \right)^{-2} \left(  k U_{k-1}\left(\frac{1+\mu}{1-\mu}\right) \right) \left(  \frac{2}{(1-\mu)^2} \right).
	\end{align*}
	Since $T_k(1) = 1$ and $U_{k-1}(1) = k$, we obtain that
	\begin{equation*}
		\lim_{ \mu \mapsto 0^+} \frac{\mu}{ 1 - p_{k}(\mu)} = \lim_{ \mu \mapsto 0^+} \frac{1}{- p_k'(\mu)} = \frac{1}{2 (2^{-2}) (k^2) (2)} = \frac{1}{k^2}.
	\end{equation*}
	This implies \eqref{ine:TG-rate-bound-chebyshev} and completes the proof. 
\end{pf}

\begin{remark}
	When $k=2$, \eqref{ine:TG-rate-bound-chebyshev} shows that the AMLI-cycle using Chebyshev polynomial converges uniformly if $\delta_{TG} < 0.75$, which is an improvement over the $\delta_{TG} < 0.5$ requirement for the W-cycle method.  Furthermore, comparing~\eqref{ine:TG-rate-bound-chebyshev} with~\eqref{ine:TG-rate-bound-simple}, it is evident that using the Chebyshev polynomials is always advantageous for any $k \geq 2$. 
\end{remark}

\begin{remark}
	According to \eqref{def:chebyshev-mu}, it is not necessary to  estimate the extreme eigenvalues on all levels.  Instead, we only need to determine $\mu$ using $\delta_{TG}$, where an upper bound is sufficient.  In practice, when the two-grid method converges uniformly, an estimation of $\delta_{TG}$ can be easily obtained on the second coarsest level, allowing $\mu$ to be determined by solving \eqref{ine:chebyshev-condition_mu}.  Alternatively, theoretical upper bound, such as those derived in \cite{brannick2013algebraic, napov2012algebraic}, can also be used. Overall, our new theory leads to implementation that could reduce the computational cost, especially for large-scale problems requiring many levels for an MG method.  
\end{remark}

\begin{remark}
	From the discussion before Theorem~\ref{thm:uniform-conv-AMLI-chebyshev}, we can see that the estimation of the smallest eigenvalue $\mu_{\ell+1}$ is essential here.
	To see this, simply using $\mu_{\ell+1} = 0$ 
	we have $ 0 \leq p_k(x) \leq p_k(0) = 1$ which implies $\mu = 0$ and $\kappa(B_{\ell}A_{\ell})$ is unbounded.  In practice, we observe that the performance of the AMLI-cycle method using the Chebyshev polynomials deteriorates if we set $\mu=0$ (see Section~\ref{sec:num}), which motivates us to seek other alternatives. 
\end{remark}

\section{AMLI-cycle using the Momentum Accelerated Polynomials}\label{sec:AMLI-Cycles-MA}
As discussed in Section~\ref{sec:AMLI-Cycle-Chebyshev}, although we introduced a new way to implement the AMLI-cycle using the Chebyshev polynomial without estimating the extreme eigenvalues on all levels, we still need to estimate $\delta_{TG}$ and then compute $\mu$ in the implementation.  In this section, we propose to use the momentum-accelerated polynomial \eqref{eqn:recurrence-relation-of-pk-MA}.  This approach eliminates the need for any estimations. Furthermore, when $k=2$, we show that the resulting momentum-accelerated AMLI-cycle is asymptotically as effective as the Chebyshev-based AMLI-cycle. 

\subsection{Analysis of AMLI-cycle without Extreme Eigenvalues}
In this subsection, we provide a new theoretical analysis of the AMLI-cycle that does not rely on extreme eigenvalues.  Our analysis is based on the following assumption about the polynomial. 
\begin{assumption}\label{asmp:bar-mu}
	For polynomial $\widetilde{p}_k(\widetilde{x})$, $\widetilde{x} \in (0, 1]$ with $\widetilde{p}_k(0) = 1$, there exists $0< \widetilde{\nu} < 1$ such that $p_k(\widetilde{\nu}) \geq 0$, $p_k(\widetilde{x}) \ \text{is decreasing on} \ (0, \widetilde{\nu}]$, and $\displaystyle p_k(\widetilde{\nu}) \geq \max_{\widetilde{x} \in [\widetilde{\nu}, 1]} [p_k(\widetilde{x})]$.
\end{assumption}
An immediate consequence of Assumption~\ref{asmp:bar-mu} is
$ \displaystyle 
\widetilde{p}_k(\nu) = \max_{\widetilde{x} \in [\nu, 1]} [\widetilde{p}_k(\widetilde{x})], \  \forall \, \nu \in (0, \widetilde{\nu})
$.

Based on Assumption~\ref{asmp:bar-mu},  we present the following condition number estimation of the AMLI-cycle, which does not require the estimation of extreme eigenvalues. 
\begin{theorem}  \label{thm:AMLI-general-p}
	Assume that the two-grid method converges uniformly with convergence rate $\delta_{TG}$ and Assumption~\ref{asmp:bar-mu} holds. The AMLI operator $B_{\ell}$, defined by AMLI-cycle (Algorithm~\ref{alg:AMLIcycle}) using polynomial $\displaystyle p_k(x) = \widetilde{p}_k(\frac{x}{L})$ with $\displaystyle L= \max\{1,  1 - \min_{\widetilde{x}\in (0, 1]} [\widetilde{p}_k(\widetilde{x})]$\}, satisfies the following condition number estimate $\kappa(B_{\ell}A_{\ell}) \leq L/\mu$
	provided that there exists $0 < \mu \leq \bar{\mu} := \min\{L \widetilde{\nu}, 1-\delta_{TG}  \}$ satisfying
	\begin{equation} \label{ine:cond-mu}
		\mu \leq (1-\delta_{TG}) \left[ 1 - p_k(\mu) \right].
	\end{equation}
\end{theorem}

\begin{pf}
	First of all, it is easy to see that $p_k(x) < 1$ for $x \in (0, L]$.  Therefore, by Theorem~\ref{thm:AMLI-SPD}, the AMLI-cycle using $p_k(x)$ is well-defined and $B_{\ell}$ are SPD for all $\ell$. 
	
	Based on ~Algorithm~\ref{alg:AMLIcycle}, we have that $B_{J} = A_{J}^{-1}$, so $\mu_J = L_J = 1$ which implies $0 < \mu \leq \mu_J = L_J \leq L$.  On level $\ell=J-1$, we have a two-grid method, i.e, $B_{J-1} = B_{J-1}^{TG}$. Thus $0< \mu_{J-1} = 1 - \delta_{TG}$ and $L_{J-1} = 1$.  By the assumptions, we have $ 0< \mu \leq \bar{\mu} \leq 1 - \delta_{TG} =  \mu_{J-1} < L_{J-1}  = 1 \leq L$.
	
	We proceed with mathematical induction. Assume $0 < \mu \leq \mu_{\ell+1} \leq L_{\ell+1} \leq L $. On level $\ell$, by \eqref{def:mu_L_recursive}, $L_{\ell} = \max\{ 1, \gamma_1^{\ell+1} \}$ and $\mu_{\ell} = (1 - \delta_{TG}) \min\{ 1, \gamma_0^{\ell+1} \}$.  Note that
	\begin{equation*}
		\gamma_1^{\ell+1}  = 1 - \min_{x \in [\mu_{\ell+1}, L_{\ell+1}]} [p_k(x)]  \leq 1 - \min_{x \in (0, L]} [p_k(x)] = 1 - \min_{\widetilde{x} \in (0,1]}[\widetilde{p}_k(\widetilde{x})] \leq L.
	\end{equation*}
	Therefore $ L_{\ell}  = \max\{ 1, \gamma_1^{\ell+1} \} \leq \max\{ 1, L \} = L$. On the other hand, for $\mu_{\ell}$, we have
	\begin{align*}
		\mu_{\ell} & = (1 - \delta_{TG}) \min\{ 1, \gamma_0^{\ell+1} \} =  (1 - \delta_{TG}) \min\{ 1,  1 - \max_{x\in [\mu_{\ell+1}, L_{\ell+1}]} [p_k(x)] \} \\
		&\geq (1 - \delta_{TG}) \left[1 - \max_{x\in [\mu,L]} [p_k(x)] \right] =  (1 - \delta_{TG}) \left[1 - \max_{\widetilde{x}\in [\frac{\mu}{L}, 1]} [\widetilde{p}_k(\widetilde{x})] \right] \\
		& = (1 - \delta_{TG}) \left[1 -  \widetilde{p}_k(\frac{\mu}{L})\right] 
		= (1 - \delta_{TG}) \left[1 -  p_k(\mu)\right] \geq \mu. 
	\end{align*}
	Here we use the fact that $0 < \frac{\mu}{L} \leq \widetilde{\nu}$ and the condition \eqref{ine:cond-mu} in the last inequality.  Thus, we have $0 < \mu \leq \mu_{\ell}\leq L_{\ell} \leq L$.  By mathematical induction, we have $0 < \mu \leq \mu_{\ell}\leq L_{\ell} \leq L$ for all $\ell$.  By Theorem~\ref{thm:AMLI-convergence}, more precisely~\eqref{ine:B_uniform_bound}, we have
	\begin{align*}
		\kappa(B_{\ell} A_{\ell}) \leq \frac{1}{1 - \delta_{TG}} \frac{\max\{ 1, \gamma_1 \}}{ \min \{ 1, \gamma_0 \}}.
	\end{align*}
	Note that 
	\begin{align*}
		\gamma_1 & = 1 - \min_{x\in [\mu, L]} [p_k(x)] \leq 1 - \min_{x\in (0,L]} [p_k(x)] =  1 - \min_{\widetilde{x}\in(0,1]} [\widetilde{p}_k(\widetilde{x})]\leq L, \\
		\gamma_0 & = 1 - \max_{x \in [\mu,L]}[p_k(x)] \geq 1 - \max_{\widetilde{x} \in [\frac{\mu}{L},1]} [\widetilde{p}_k(\widetilde{x})] = 1 - \widetilde{p}_k(\frac{\mu}{L}) = 1 - p_k(\mu).
	\end{align*}
	Thus, 
	$\kappa(B_{\ell} A_{\ell}) \leq \frac{L}{ (1-\delta_{TG}) [1 - p_k(\mu)] } \leq \frac{L}{\mu}$,
	which completes the proof. 
\end{pf}

We emphasize that the uniformly bounded condition number result of the AMLI-cycle, derived in~Theorem~\ref{thm:AMLI-general-p}, does not depend on estimating extreme eigenvalues. Once the polynomial $\widetilde{p}_k(\widetilde{x})$ is specified on $[0,1]$, we can determine $L$, and from there, the polynomial $p_k(x)$ is fully defined. The uniform condition number estimation result holds as long as the conditions outlined in Theorem~\ref{thm:AMLI-general-p} are satisfied.  

\begin{remark}
	Theorem~\ref{thm:AMLI-general-p} also provides insight into why the AMLI-cycle using Chebyshev polynomial might not work if we simply use $\mu_{\ell} = 0$.  In this case, we have $L=1$ and, thus, $p_k(x) = \widetilde{p}_k(\widetilde{x})$.   If $\mu_{\ell}=0$, then Assumption~\ref{asmp:bar-mu} does not hold because $p_k(x) = 1$ for some $x \in (0, 1]$, meaning that $\widetilde{\nu}$ does not exist.  
\end{remark}

\subsection{Analysis of the Momentum-accelerated AMLI-cycle}
For the polynomials derived from the momentum acceleration methods, i.e., \eqref{eqn:recurrence-relation-of-pk-MA},  we simply set $\mu_{\ell+1} = 0$ and $L_{\ell+1} = L$ in our derivation to avoid the need for estimating extreme values. Here, $L$ is a predetermined paramter (see Theorem~\ref{thm:AMLI-general-p}).  From Table~\ref{tab:choices-MA},  for both NA and sAA(1) methods, we have $\beta = 1$ when $\mu_{\ell+1}=0$. Thus, we use $\beta = 1$ in the MA polynomial~\eqref{eqn:recurrence-relation-of-pk-MA}. Regarding the choice of $\alpha$, when $\mu_{\ell+1} = 0$, from Table~\ref{tab:choices-MA}, we have $\alpha = \frac{1}{L}$ for the NA method and  $\alpha = \frac{4}{3L}$ for the sAA(1) method.  Thus, we choose $\alpha = \frac{a}{L}$ and discuss different choices of $a$.  These choices yield
$p_{k+1}(x) = 2 \left(1 - a \frac{x}{L}\right) p_k(x) - \left(1-a \frac{x}{L} \right) p_{k-1}(x)$.
Therefore, it is easy to see that the corresponding $\widetilde{p}_k(\widetilde{x})$ is defined as follows,
\begin{equation} \label{def:tilde-pk}
	\widetilde{p}_{k+1}(\widetilde{x}) = 2 \left(1 - a \widetilde{x}\right) \widetilde{p}_k(\widetilde{x}) - \left(1-a \widetilde{x} \right) \widetilde{p}_{k-1}(\widetilde{x}).
\end{equation}
In addition, we set $\widetilde{p}_0(\widetilde{x}) = 1$ and $\widetilde{p}_1(\widetilde{x}) = 1 - \widetilde{x}$, which are the same as the first two polynomials of the simple polynomial and the scaled and shifted Chebyshev polynomials.  

Next, we apply Theorem~\ref{thm:AMLI-general-p} to the AMLI-cycle using the MA polynomial \eqref{def:tilde-pk} and discuss different choices of $k$ and $a$.

\subsubsection{Case: $k=2$} When $k=2$, by direct calculation, $\widetilde{p}_2(\widetilde{x})$ is defined as
\begin{equation}\label{def:MA-tilde-p2}
	\widetilde{p}_2(\widetilde{x})  = (1 - a\widetilde{x}  ) \left( 1 - 2 \widetilde{x}  \right).
\end{equation}
First of all, to make sure $\widetilde{p}_2(\widetilde{x})<1$ for $\widetilde{x} \in (0,1]$, we need $0 < a < 2$.   Next lemma gives the minimum value of $\widetilde{p}_2(\widetilde{x})$ on interval $(0, 1]$. 
\begin{lemma}\label{lem:p2-min}
	For the MA polynomial of degree $2$ defined in \eqref{def:MA-tilde-p2}, we have
	\begin{equation*}
		\min_{\widetilde{x} \in (0,1]} [ \widetilde{p}_2(\widetilde{x}) ] =
		\begin{cases}
			\widetilde{p}_2(1) = a-1, &\quad 0 < a < \frac{2}{3}, \\
			\widetilde{p}_2(\frac{2+a}{4a}) = -\frac{(2-a)^2}{8a}, 
			&\quad \frac{2}{3} \leq a < 2.
		\end{cases}
	\end{equation*} 
\end{lemma}
\begin{pf}
	The minimal value of $p_k(x)$ can be obtained by direct calculation. 
\end{pf}

To apply Theorem~\ref{thm:AMLI-general-p}, we need to verify Assumption~\ref{asmp:bar-mu}, which leads to the next lemma. 
\begin{lemma} \label{lem:p2-asmp}
	Assumption~\ref{asmp:bar-mu} holds for $\widetilde{p}_2(\widetilde{x})$ with $\displaystyle 0 < a < 2$.
\end{lemma}
\begin{pf}
	It is easy to see that $\widetilde{p}_2(\widetilde{x})$ is decreasing on $(0, \frac{1}{2}]$. Thus, we can simply choose $0< \widetilde{\nu} \leq \frac{1}{2}$ such that $\widetilde{p}_2(\widetilde{\nu}) = \max_{\widetilde{x} \in [\frac{1}{2},1]} [ \widetilde{p}_2(\widetilde{x})]$.  This choice of $\widetilde{\nu}$ satisfies Assumption~\ref{asmp:bar-mu}. 
\end{pf}

Based on Lemma~\ref{lem:p2-min}, Lemma~\ref{lem:p2-asmp}, and Theorem~\ref{thm:AMLI-general-p}, we have the following corollary regarding the condition number of the AMLI-cycle using the MA polynomial when $k=2$.
\begin{corollary}\label{coro:AMLI-p2}
	Assuming that the two-grid method converges uniformly with convergence rate $\delta_{TG}$. The AMLI operator $B_{\ell}$, defined by the AMLI-cycle (Algorithm~\ref{alg:AMLIcycle}) using the MA polynomial $p_2(x):= \widetilde{p}_2(\frac{x}{L})$ with 
	\begin{equation*}
		L= 
		\begin{cases}
			2- a, & \quad 0 < a < \frac{2}{3}, \\
			\frac{(2+a)^2}{8a}, & \quad \frac{2}{3} \leq a < 2, 
		\end{cases}
	\end{equation*}
	satisfies the condition number estimate $\kappa(B_{\ell}A_{\ell}) \leq \frac{L}{\mu}$ provided there exists $0 < \mu \leq \bar{\mu} := \min\{L \widetilde{\nu}, 1-\delta_{TG}  \}$ such that, 
	\begin{equation} \label{ine:cond-mu-p2}
		\mu \leq (1-\delta_{TG}) \left[ 1 - p_2(\mu) \right].
	\end{equation}
	Moreover, a sufficient condition for the existence of $\mu$ in \eqref{ine:cond-mu-p2} is
	\begin{equation*} 
		\delta_{TG} < 
		\begin{cases}
			1 - \frac{2-a}{2+a}, & \quad 0 < a < \frac{2}{3}, \\
			1 - \frac{2+a}{8a}, & \quad \frac{2}{3} \leq a < 2.
		\end{cases}
	\end{equation*}
\end{corollary}
\begin{pf}
	We only discuss the case $\frac{2}{3} \leq a < 2$. The result for $0 < a < \frac{2}{3}$ can be obtained by the same argument.  By Lemma~\ref{lem:p2-min} and Theorem~\ref{thm:AMLI-general-p}, we have $L = 1 - \min_{\widetilde{x} \in (0,1]} [\widetilde{p}_2(\widetilde{x})] = 1 + \frac{(2-a)^2}{8a} = \frac{(2+a)^2}{8a}$.  Thus, 
	\begin{equation*}
		p_2(x) = \left(1 - a \frac{x}{L}\right)\left(1-2\frac{x}{L} \right) = \left( 1 -  \frac{8a^2x}{(2+a)^2} \right) \left( 1 -  \frac{16ax}{(2+a)^2}  \right).
	\end{equation*}	
	The condition \eqref{ine:cond-mu-p2} becomes,
	\begin{equation*}
		\mu \leq \left[ 1- \left(1 -  \frac{8a^2\mu}{(2+a)^2} \right) \left( 1 - \frac{ 16a\mu}{(2+a)^2} \right) \right] (1 - \delta_{TG}).
	\end{equation*} 
	This implies the existence of $\mu$ if 
	\begin{equation*}
		\delta_{TG} \leq 1 - \frac{1}{\frac{8a}{2+a} - \frac{128a^3\mu}{(2+a)^4} } < 1 - \frac{2+a}{8a},
	\end{equation*}
	which completes the proof. 
\end{pf}

\begin{remark} \label{remark:p2-best}
	From \eqref{ine:cond-mu-p2}, we have, for $\frac{2}{3} \leq a  < 2$,
	\begin{equation*}
		\delta_{TG} < 1 - \frac{2+a}{8a} < \frac{3}{4},
	\end{equation*}
	where the upper bound is obtained as  $a \mapsto2$.  This implies that the AMLI-cycle using the MA polynomial is uniform preconditioner if $\delta_{TG} < 0.75$ as $a \mapsto 2$, which matches the convergence behavior of the AMLI-cycle using the Chebyshev polynomials.  In the practical implementation, we set $a$ close to $2$, and there is no need to estimate $\delta_{TG}$.   Additionally, if $a =  \frac{2}{3}$, we obtain $\delta_{TG} < \frac{1}{2}$, which is the same condition as the W-cycle method. This means we need to consider the range $\frac{2}{3} < a < 2$ to achieve improvement.
\end{remark}

\subsubsection{Case:  $k=3$} When $k=3$, based on~\eqref{def:tilde-pk}, the MA polynomial is 
\begin{equation} \label{def:MA-tilde-p3}  
	\widetilde{p}_3(\widetilde{x}) = 2 \left(  1-  a \widetilde{x} \right)^2(1 -  2 \widetilde{x}) - \left( 1 - a \widetilde{x} \right) \left( 1 - \widetilde{x} \right).
\end{equation}
By direct calculation, we still have $\widetilde{p}_3(\widetilde{x}) < 1$, $\widetilde{x} \in (0,1]$, for $0 < a < 2$.  Next  lemma discusses its minimum value. 
\begin{lemma}
	For the MA polynomial of degree $3$ defined by~\eqref{def:MA-tilde-p3}, we have
	\begin{equation*}
		\min_{\widetilde{x} \in (0,1]} [\widetilde{p}_3(\widetilde{x})] =
		\begin{cases}
			\widetilde{p}_3(\frac{2a+7 - \sqrt{4(1-a)^2 + 9}}{12a}),  &\quad 	0 < a < \frac{9 + 2  \sqrt{22}}{14}, \\
			\widetilde{p}_3(1) = -2(a-1)^2, &\quad  \frac{9 + 2  \sqrt{22}}{14} \leq a < 2.
		\end{cases}
	\end{equation*}
\end{lemma}

Next lemma verifies Assumption~\ref{asmp:bar-mu} for $\widetilde{p}_3(\widetilde{x})$.
\begin{lemma} \label{lem:p3-asmp}
	Assumption~\ref{asmp:bar-mu} holds for $\widetilde{p}_3(\widetilde{x})$ with $\displaystyle 0 < a < 2$.
\end{lemma}
\begin{pf}
	It is easy to see that $\widetilde{p}_2(\widetilde{x})$ is decreasing on $(0, \frac{2a+3-\sqrt{(2a-1)^2+8}}{8a}]$ for $0 < a < 2$ (note $\frac{2a+3-\sqrt{(2a-1)^2+8}}{8a}$ is the smallest root located in $(0,1]$).  Then the existence of $\widetilde{\nu}$ follows from the same argument of the proof of Lemma~\ref{lem:p2-asmp}.
\end{pf}

Now we are ready to present the condition number of the AMLI-cycle using the MA polynomial $\widetilde{p}_3(\widetilde{x})$ in the next corolloary.

\begin{corollary}\label{coro:AMLI-p3}
	Assuming that the two-grid method converges uniformly with convergence rate $\delta_{TG}$. The AMLI operator $B_{\ell}$, defined by the AMLI-cycle (Algorithm~\ref{alg:AMLIcycle}) using the MA polynomial $p_3(x):= \widetilde{p}_3(\frac{x}{L})$ with 
	\begin{equation*}
		L= 
		\begin{cases}
			1 - \widetilde{p}_3(\frac{2a+7 - \sqrt{4(1-a)^2 + 9}}{12a}),  &\quad 	0 < a < \frac{9 + 2  \sqrt{22}}{14},  \\
			1+ 2(a-1)^2, & \quad \frac{9 + 2  \sqrt{22}}{14} \leq a < 2, 
		\end{cases}
	\end{equation*}
	satisfies the condition number estimate $\kappa(B_{\ell}A_{\ell}) \leq \frac{L}{\mu}$ provided there exist $0 < \mu \leq \bar{\mu} := \min\{L \widetilde{\nu}, 1-\delta_{TG}  \}$ such that, 
	\begin{equation} \label{ine:cond-mu-p3}
		\mu \leq (1-\delta_{TG}) \left[ 1 - p_3(\mu) \right].
	\end{equation}
	Moreover, a sufficient condition for the existence of $\mu$ in \eqref{ine:cond-mu-p3} is
	\begin{equation}\label{ine:p3-mu-TG-condition}
		\delta_{TG} < 
		\begin{cases}
			1 - \frac{1 - \widetilde{p}_3(\frac{2a+7 - \sqrt{4(1-a)^2 + 9}}{12a})}{3a+3}, & \quad 0 < a < \frac{9 + 2  \sqrt{22}}{14}, \\
			1 - \frac{1+2(a-1)^2}{3a+3}, & \quad \frac{9 + 2  \sqrt{22}}{14} \leq a < 2. 
		\end{cases}
	\end{equation}
\end{corollary}

\begin{pf}
	The proof of the condition number and \eqref{ine:cond-mu-p3} follows from the same argument of the proof of Corollary~\ref{coro:AMLI-p2}.   To show \eqref{ine:p3-mu-TG-condition}, note that we can rewrite \eqref{ine:cond-mu-p3} as 
	\begin{equation*}
		\widetilde{\mu} \leq \frac{1-\delta_{TG}}{L} [1 - \widetilde{p}_3(\widetilde{\mu})],
	\end{equation*}
	where $\widetilde{\mu} := \frac{\mu}{L}$.  Thus, we have 
	\begin{align*}
		\delta_{TG} &\leq 1 - \frac{L \widetilde{\mu}}{ 1- \widetilde{p}_3(\widetilde{\mu})}= 1 - \frac{L }{(3a+3) - a(2a+7) \widetilde{\mu} + 4a^2 \widetilde{\mu}^2}  \leq 1 - \frac{L}{3a+3}.
	\end{align*}
	This completes the proof. 
\end{pf}

\begin{remark}
	Let us examine the sufficient condition \eqref{ine:p3-mu-TG-condition} to find the optimal upper bound of $\delta_{TG}$.  It turns out, unlike the case for $k=2$, the best upper bound of $\delta_{TG}$ is obtained when $a = \frac{9 + 2 \sqrt{22}}{14}$ and 
	$\delta_{TG} <  \frac{8 + 2\sqrt{22}}{21} \approx 0.827659$. 
	Compared to the AMLI-cycle using the Chebyshev polynomials, which requires $\delta_{TG} < \frac{8}{9} \approx 0.888889$ when $k=3$,  the AMLI-cycle using the MA polynomial is nearly optimal and very close to this bound when $a = \frac{9 + 2 \sqrt{22}}{14}$.  
\end{remark}

\subsection{Case: General $k$}
Now, let us consider the general case for $k$. Unfortunately, starting from $k=4$, $\widetilde{p}_k(\widetilde{x})$ might be bigger than $1$ for $0 < a <2$.  In what follows, we drive the range of $\alpha$ such that $|\widetilde{p}_k(\widetilde{x})| < 1$, which is stronger than the minimal requirement, as it also provides a lower bound of $\widetilde{p}_k(\widetilde{x})$, i.e., $\widetilde{p}_k(\widetilde{x}) > -1$.   The result is summarized in the following lemma. 
\begin{lemma}\label{lem:pk<1}
	When $k\geq 4$ and $\frac{1}{2} \leq a \leq \frac{4}{3}$, the MA polynomial $\widetilde{p}_k(\widetilde{x})$ \eqref{def:tilde-pk} satisfies 
	\begin{equation*}
		\max_{\widetilde{x} \in (0,1]}|\widetilde{p}_k(\widetilde{x})| < 1.
	\end{equation*} 
\end{lemma}
The proof is based on rewriting the three-term recurrence relationship \eqref{def:tilde-pk} into a two-term one using a $2 \times 2$ matrix form.  The details are in \ref{appendix:proof-pk<1}.

Next, we estimate the condition number of the AMLI-cycle when $\frac{1}{2} \leq a \leq \frac{4}{3}$, under the assumption that the two-grid method converges uniformly with convergence rate $\delta_{TG}$.  To proceed, we first need the following lemma to verify Assumption~\ref{asmp:bar-mu}.

\begin{lemma}\label{lem:pk-asmp}
	For $\frac{1}{2} \leq a \leq \frac{4}{3}$, Assumption~\ref{asmp:bar-mu} holds for $\widetilde{p}_k(\widetilde{x})$ with $k \geq 4$.
\end{lemma}
\begin{pf}
	Since $\widetilde{p}_k(0) = 1$ and $\max_{\widetilde{x} \in (0,1]} |\widetilde{p}_k(\widetilde{x})| < 1$ by Lemma~\ref{lem:pk<1},  as a polynomial, if the smallest positive root $r_0$ of $\widetilde{p}_k(\widetilde{x})$ is located inside $(0,1]$, then $\widetilde{p}_k(\widetilde{x})$ is decreasing on $(0,r_0]$ and the existence of~$\widetilde{\nu}$ can be obtained by the same argument of Lemma~\ref{lem:p2-asmp}.  Otherwise, there is no root located inside $(0, 1]$ and $\widetilde{p}_k(\widetilde{x})$ is decreasing on $(0,1]$.  In this case, we can choose any $0 < \widetilde{\nu} < 1$, and this choice satisfies Assumption~\ref{asmp:bar-mu}.
\end{pf}

Based on Lemma~\ref{lem:pk-asmp},  we can apply Theorem~\ref{thm:AMLI-general-p} and have the following corollary showing the uniformly bounded condition number. 

\begin{corollary}\label{coro:AMLI-1/2<=alpha<=4/3}
	Assume that the two-grid method converges uniformly with convergence rate $\delta_{TG}$. The AMLI operator $B_{\ell}$, defined by the AMLI-cycle (Algorithm~\ref{alg:AMLIcycle}) using the MA polynomial~\eqref{def:tilde-pk} $p_k(x) := \widetilde{p}_k(\frac{x}{L})$ with $L=2$ and $\frac{1}{2} \leq a \leq \frac{4}{3} $ satisfies the condition number estimate $	\kappa(B_{\ell}A_{\ell}) \leq \frac{2}{\mu}$, provided  there exists $0 < \mu \leq \bar{\mu}:=\min\{ L \widetilde{\nu}, 1-\delta_{TG} \}$ satisfying \eqref{ine:cond-mu}.
\end{corollary}
\begin{pf}
	The result follows directly from Theorem~\ref{thm:AMLI-general-p} with $L=2$.
\end{pf}

\begin{remark}
	In Corollary~\ref{coro:AMLI-1/2<=alpha<=4/3}, the upper bound $L$ is relatively loose, leading to a suboptimal conclusion.   To illustrate this, condition \eqref{ine:cond-mu} implies the following requirement on the two-grid method convergence rate
	$\delta_{TG} < 1 - \frac{1}{\frac{k(k-1)}{4}a + \frac{k}{2}}$.
	When $k=2$ and $\frac{1}{2} \leq a \leq \frac{4}{3}$, we have
	$\delta_{TG} < \frac{2}{5}$.
	Comparing with the discussion in Remark~\ref{remark:p2-best}, we can see the result is suboptimal.  In practice, one can compute $\widetilde{p}_{\min} = \min_{\widetilde{x} \in (0,1]} \widetilde{p}_k(\widetilde{x}) $ and use $L = 1 - \widetilde{p}_{\min}$ for a given $k$.

\end{remark}

\subsection{Practical Implementations}
In this subsection, we present our practical implementation of the AMLI-cycle using the MA polynomial \eqref{def:tilde-pk}.  The general idea is, for small $k$, e.g., $k=2$ and $k=3$, we use the accurate upper bound $L$ provided by Corollary~\ref{coro:AMLI-p2} and Corollary~\ref{coro:AMLI-p3} and, consequently, a good choice of $a$ based on the theory.  For large $k$, we follow Corollary~\ref{coro:AMLI-1/2<=alpha<=4/3} and use $L = 2$ and $a = \frac{4}{3}$.  Although they might be suboptimal theoretically, their practical performance is nearly optimal according to our numerical experiments.

\begin{algorithm}
	\caption{Practical AMLI-cycle MG:~$B_\ell \mathbf{b}$} \label{alg:AMLIcycle-MA}
	\begin{algorithmic}[1]
		\IF{$\ell == J$}
		\STATE $\mathbf{x_\ell} = A_{\ell}^{-1} \mathbf{b}$
		\ELSE
		\STATE $\textbf{Presmoothing}$: $\mathbf{x}_{\ell} \leftarrow  M_{\ell}\mathbf{b}$
		\STATE $\textbf{Restriction}$: $\mathbf{r}_{\ell+1} \gets R_{\ell} (\mathbf{b} - A_{\ell} \mathbf{x}_{\ell})$
		\STATE $\textbf{Coarse-grid\ correction}$: 
		\IF{$\ell = J-1$}
		\STATE $\mathbf{e}_{\ell+1} \gets B_{\ell+1} \mathbf{r}_{\ell+1} $
		\ELSE 
		\STATE Set $L$ and $a$ as follows, 
		\vskip -20pt
		\begin{equation*}
			\begin{cases}
				a \gets 1.9, & \text{if\ } k=2 \\
				a \gets \frac{9 + 2 \sqrt{22}}{14},  &  \text{if\ } k=3   \\
				a \gets 4/3.   & \text{otherwise}
			\end{cases}
			\quad 
			\text{and}
			\quad
			\begin{cases}
				L \gets 1, & \text{if\ } k=1 \\
				L \gets \frac{(2+a)^2}{8a}, & \text{if\ } k=2 \\
				L \gets 1 + 2(a-1)^2, & \text{if\ } k=3 \\
				L \gets  2.  & \text{otherwise}
			\end{cases}
		\end{equation*}
		\vskip -6pt
		\STATE $\mathbf{e}_{\ell+1}^0 \gets \mathbf{0}$ 
		\STATE $\mathbf{e}_{\ell+1}^1 \gets \mathbf{e}_{\ell+1}^0 + \frac{1}{L} B_{\ell+1} (\mathbf{r}_{\ell+1} - A_{\ell+1} \mathbf{e}_{\ell+1}^0 )$ 
		\FOR{$i = 2, \cdots, k$}
		\STATE $\mathbf{e}_{\ell+1}^i     \leftarrow 2\left[\mathbf{e}_{\ell+1}^{i-1}+\frac{a}{L} B_{\ell+1}(\mathbf{r}_{\ell+1}-A\mathbf{e}_{\ell+1}^{i-1}) \right]- \left[ \mathbf{e}_{\ell+1}^{i-2}+\frac{a}{L} B_{\ell+1}(\mathbf{r}_{\ell+1}-A\mathbf{e}_{\ell+1}^{i-2}) \right]$
		\ENDFOR
		\STATE $\mathbf{e}_{\ell+1} \gets \mathbf{e}_{\ell+1}^k$
		\ENDIF
		\STATE $\textbf{Prolongation}$: $\mathbf{x}_{\ell} \gets \mathbf{x}_{\ell} + P_{\ell} \mathbf{e}_{\ell+1}$
		\STATE $\textbf{Postsmoothing}$: $\mathbf{x}_{\ell} \leftarrow \mathbf{x}_{\ell} + M_{\ell}^{T} (\mathbf{b} - A_\ell\mathbf{x}_{\ell})$
		\ENDIF
		\STATE $B_{\ell} \mathbf{b} \gets \mathbf{x}_{\ell}$
	\end{algorithmic}
\end{algorithm}

\begin{remark}
	When $k=2$, theory suggests us to use $a \mapsto 2$ but not $a=2$. Thus, in our implementation, we suggest $a = 1.9$.  
\end{remark}

\section{Numerical Results}\label{sec:num}
In this section, we present numerical experiments to demonstrate the efficiency of the proposed practical implementation of the AMLI-cycle using the scaled and shifted Chebyshev polynomial~\eqref{def:chebyshev}, as discussed in Section~\ref{sec:AMLI-Cycle-Chebyshev}, and the AMLI-cycle using the momentum acceleration polynomial~\eqref{eqn:recurrence-relation-of-pk-MA}, as discussed in Section~\ref{sec:AMLI-Cycles-MA}.  For brevity, we refer to these methods as C-AMLI-cycle and M-AMLI-cycle, respectively.  

In our numerical experiments, we use zero right-hand side, i.e., $\mathbf{b} = \mathbf{0}$, and a random initial guess. We use the unsmoothed aggregation AMG~(UA-AMG) method~\cite{Kim03} in all the experiments since it is well-known that the V-cycle UA-AMG method does not converge uniformly in general and more sophisticated cycles are needed. Our implementation of UA-AMG utilizes the maximal independent set (MIS) algorithm to construct the aggregations and define the coarse-level matrices. For smoothing, we use Gauss-Seidel (GS) smoother ($1$ step forward GS for pre-smoothing and~$1$ step backward GS for post-smoothing).  Additionally, UA-AMG is used as a preconditioner within the preconditioned conjugate gradient (PCG) method (for K-cycle, generalized PCG with variable preconditoner is used). The stopping criterion is a relative residual less than or equal to~$10^{-6}$.  All the numerical experiments were conducted on a MacBook Pro with an Apple M1 Max CPU and 64 GB of RAM. 

\subsection{Standard Poisson Equation}
We start with the standard Poisson problem as our first example to demonstrate the effectiveness of the proposed methods.
\begin{example}\label{exp:poisson}
	Consider the model problem on~$\Omega = [0, 1] \times [0, 1]$.
	\begin{equation*}
		\begin{array}{rcl}
			- \Delta \mathbf{u} &=& f, \ \  \ \text{in} \ \ \Omega,\\
			\mathbf{u} &=& 0, \ \ \ \text{on} \ \ \partial \Omega.
		\end{array}
	\end{equation*}
\end{example}

We use the standard linear finite-element method on a uniform triangulation of $\Omega$ for the discretization.  We choose $f=0$ here so that the right-hand side for the linear system is zero.  As a result, the exact solution is zero as well.  Since the problem is isotropic, coarsening is performed directly using an MIS-based aggregation scheme.  For the mesh of size $h = 1/2048$, the minimal coarsening ratio ( the smallest ratio of consecutive matrix sizes across levels) is approximately $6$.   Based on the computational complexity of the AMLI-cycle (see \cite{Vassilevski08}), we can use $k \leq 6$ to maintain optimal complexity.  In our experiment, we use $k \leq 5$.  

In~Table~\ref{tab:Exa.Poissontwogrid and cycles}, we present the number of iterations of PCG with different types of MG cycles as preconditioners for different mesh sizes.  The results highlight several trends. Although the two-grid method achieves uniform convergence, the V-cycle (i.e., $k$V-cycle with~$k=1$) exhibits performance degradation, as expected, due to the UA-AMG approach.  Since the convergence factor of the two-grid method is greater than $0.5$, W-cycle (i.e., $k$V-cycle with $k = 2$) also fails to achieve uniform convergence, resulting in an increasing number of iterations as~$h$ decreases.  For higher values of~$k$, the PCG iteration count continues to grow slightly for~$k=3$ and~$k=4$, but the results suggest that uniform convergence is attained for~$k=5$. We include the K-cycle in our comparison due to its known optimality among AMLI-type MG cycles \cite{Hu13}. As shown in~Table~\ref{tab:Exa.Poissontwogrid and cycles}, the K-cycle achieves uniform convergence starting at~$k=2$ and performs as good as the two-grid method when $k\ge 3$. Next, we evaluate our implementation of the C-AMLI-cycle, which requires only an estimation of the two-grid convergence $\delta_{TG}$.  As discussed in Section~\ref{sec:AMLI-Cycle-Chebyshev}, an estimation on the coarse mesh is sufficient.  For this example, we compute the two-grid convergence rate on a mesh with $h=1/128$, obtaining approximately $0.7233$.  Accordingly, we set $\delta_{TG} = 0.725$.  As shown in~Table~\ref{tab:Exa.Poissontwogrid and cycles}, the resulting C-AMLI-cycle performs effectively and comparable with the K-cycle for all tested values of $k$ and $h$.   However, if we do not estimate $\delta_{TG}$ and use $\delta_{TG} = 1$ (i.e., $\mu=0$ based on \eqref{ine:chebyshev-condition_mu}), although the number of iterations seems to be nearly uniform for $k = 2$, the performance of the C-AMLI-cycle deteriorates for $k \geq 3$.  This observation is consistent with the theoretical insights discussed in Section~\ref{sec:AMLI-Cycle-Chebyshev}.  Finally, for the M-AMLI-cycle, we set $L$  and $a$ as presented in~Algorithm~\ref{alg:AMLIcycle-MA}.   When $k=2$, the number of iterations is slightly higher than the K-cycle but comparable with the C-AMLI-cycle.  Furthermore, for  $k \geq 3$, the M-AMLI-cycle achieves uniform convergence, with performance matching that of the two-grid and the K-cycle. These results validate the theoretical predictions presented in~Theorem~\ref{thm:AMLI-general-p}.

\begin{table}[H]
	\begin{center}
		\caption{Number of PCG iterations for Example~\ref{exp:poisson}.}
		\label{tab:Exa.Poissontwogrid and cycles}
		\renewcommand{\arraystretch}{1.0}
		\begin{tabular}{cccccc}
			\hline \hline
			&  $h=1/128$ & $h=1/256$ & $h=1/512$ & $h=1/1024$ & $h=1/2048$ \\ \cline{2-6}
			&\multicolumn{5}{c}{Two-grid}  \\  \cline{2-6}
			& 11 & 11 & 11 & 11 & 11\\ \hline
			&\multicolumn{5}{c}{ $k$V-cycle  }  \\ \cline{2-6}
			$k=1$     & 25 & 37 & 44 & 61& 82\\
			$k=2$     & 17  & 19 & 22 & 24 & 26\\
			$k=3$     & 14  & 15 & 16 & 16 & 16\\
			$k=4$     & 12  & 13 & 13 & 14 & 14\\
			$k=5$     & 12  & 12 & 12 & 13 & 12\\
			\hline
			&\multicolumn{5}{c}{K-cycle }  \\ \cline{2-6}
			$k=2$     & 12 & 12 & 12 & 12 & 12\\
			$k=3$     & 11 & 11 & 11 & 11 & 11\\
			$k=4$     & 11 & 11 & 11 & 11 & 11\\
			$k=5$     & 11 &  11 & 11 & 11 & 11\\
			\hline
			&\multicolumn{5}{c}{C-AMLI-cycle ($\delta_{TG} = 0.725$)}  \\ \cline{2-6}
			$k=2$     & 12 & 13 & 13 & 14 & 13\\
			$k=3$     & 11 & 11  & 12 & 12 & 12\\
			$k=4$     & 11 & 11  & 11 & 11 & 11 \\
			$k=5$     & 11 & 11  & 11 & 11 & 11 \\
			\hline
			&\multicolumn{5}{c}{C-AMLI-cycle ($\delta_{TG}=1$) }  \\ \cline{2-6}
			$k=2$     & 12 & 13  & 13 & 13 & 13\\
			$k=3$     & 15 & 15  & 15 & 15 & 15 \\
			$k=4$     & 23 & 23 & 22 & 23 & 22 \\
			$k=5$     & 30 & 31  & 33 & 42 & 42\\
			\hline
			&\multicolumn{5}{c}{M-AMLI-cycle}  \\ \cline{2-6}
			$k=2$     & 12 & 13 & 13 & 14 & 13\\
			$k=3$     & 11 & 11 & 11  & 11 & 11 \\
			$k=4$     & 10 & 11 & 11 & 11 & 11\\
			$k=5$     & 10 & 10 & 10 & 10 & 10\\
			\hline  \hline
		\end{tabular}
	\end{center}
\end{table}

In Table~\ref{tab:Exa.Poissontwogrid and cycles-CPUtime}, we present the CPU times for the K-cycle, C-AMLI-cycle (with $\delta_{TG} = 0.725$), and M-AMLI-cycle, as these methods exhibit stable iteration counts. As expected, both the C-AMLI-cycle and M-AMLI-cycle outperform the K-cycle in terms of computational efficiency. This is primarily due to the nonlinear nature of the K-cycle, which requires the use of generalized PCG, introducing additional computational overhead. On average, the C-AMLI-cycle and M-AMLI-cycle are approximately twice as fast as the K-cycle, with comparable CPU times between them. Given that the M-AMLI-cycle does not require estimating the two-grid convergence rate, our numerical results suggest that it is a better choice in practice.

\begin{table}[H]
	\begin{center}
		\caption{CPU time for Example~\ref{exp:poisson}.}
		\label{tab:Exa.Poissontwogrid and cycles-CPUtime}
		\renewcommand{\arraystretch}{1.0}
		\begin{tabular}{cccccc}
			\hline \hline
			&  $h=1/128$ & $h=1/256$ & $h=1/512$ & $h=1/1024$ & $h=1/2048$ \\ \cline{2-6}
			&\multicolumn{5}{c}{K-cycle }  \\ \cline{2-6}
			$k=2$     & 0.047 & 0.128 & 0.539 & 1.955 & 7.238 \\
			$k=3$     & 0.051 & 0.147 & 0.567 & 2.052 & 7.803 \\
			$k=4$     & 0.057 & 0.212 & 0.735 & 2.466 & 10.229\\
			$k=5$     & 0.071 & 0.253 & 0.797 & 3.213 &  13.495\\
			\hline
			&\multicolumn{5}{c}{C-AMLI-cycle ($\delta_{TG} = 0.725$)}  \\ \cline{2-6}
			$k=2$     & 0.016 & 0.071 & 0.289 & 1.075 & 4.117 \\
			$k=3$     & 0.017 & 0.074  & 0.299 & 1.169 & 4.417  \\
			$k=4$     & 0.021 & 0.100  & 0.356 & 1.339 & 5.146 \\
			$k=5$     & 0.030 & 0.121  & 0.443 & 1.693 & 6.901 \\
			\hline
			&\multicolumn{5}{c}{M-AMLI-cycle}  \\ \cline{2-6}
			$k=2$     & 0.015 & 0.078 & 0.269 & 1.086 & 3.797 \\
			$k=3$     & 0.017 & 0.070 & 0.278 & 1.029 & 4.001 \\
			$k=4$     & 0.019 & 0.092 & 0.352 & 1.300 & 5.126 \\
			$k=5$     & 0.022 & 0.105 & 0.388 & 1.519 & 6.335 \\
			\hline  \hline
		\end{tabular}
	\end{center}
\end{table}

\subsection{Anisotropic Diffusion Problem}
The second example we consider is an anisotropic diffusion problem.
\begin{example}\label{exp:aniso-poisson}
	Let $\Omega = [0, 1] \times [0,1]$, consider
	\begin{equation*}
		\begin{array}{rcl}
			- \partial _{\mathbf{x}\mathbf{x}}\mathbf{u}- 10^{-3} \partial _{\mathbf{y}\mathbf{y}}\mathbf{u} &=& f, \ \ \ \text{in} \ \ \Omega,\\
			\mathbf{u} &=& 0, \ \ \ \text{on} \ \ \partial \Omega.
		\end{array}
	\end{equation*}
\end{example}

We again apply the standard linear finite-element method on a uniform triangulation of $\Omega$ for discretization. Again, we set $f=0$, which leads to a zero exact solution. Since the problem is anisotropic, in the coarsening step, we first drop some entries of the linear system based on the so-called strength connection (see \cite{Vassilevski08,Xu17} and we use $0.25$ for the strength connection parameter) and then perform the MIS-based aggregation scheme. On the mesh of size $h = 1/2048$, the minimal coarsening ratio is about $3$, which suggests $k \leq 3$. However,  since the average coarsening ratio is $5.6054$, we test $k \leq 4$ for this example to maintain nearly optimal complexity. In addition, we compute the two-grid convergence rate on $h=1/128$, which is about $0.7130$.  Consequently, we use $\delta_{TG} = 0.715$ for this example. 

\begin{table}[H]
	\begin{center}
		\caption{Number of PCG iterations for Example~\ref{exp:aniso-poisson}.} 
		\label{tab:Exa.aniso-poissontwogrid and cycles}
		\renewcommand{\arraystretch}{1.0}
		\begin{tabular}{cccccc}
			\hline \hline
			&  $h=1/128$ & $h=1/256$ & $h=1/512$ & $h=1/1024$ & $h=1/2048$ \\ \hline
			&\multicolumn{5}{c}{Two-grid}  \\  \cline{2-6}
			& 10 & 10 & 10 & 10 & 10 \\ 
			\hline
			&\multicolumn{5}{c}{ $k$V-cycle  }  \\ \cline{2-6}
			$k=1$     & 25 & 35 & 44 & 53 & 64 \\
			$k=2$     & 16 & 19 & 22 & 24 & 26 \\
			$k=3$     & 13 & 15 & 16 & 18 & 18 \\
			$k=4$     & 12 & 12 & 13 & 14 & 15 \\
			\hline
			&\multicolumn{5}{c}{K-cycle }  \\ \cline{2-6}
			$k=2$     & 12 & 12 & 12 & 12 & 12 \\
			$k=3$     & 11 & 11 & 11 & 11 & 11 \\
			$k=4$     & 11 & 10 & 11 & 11 & 11 \\
			\hline
			&\multicolumn{5}{c}{C-AMLI-cycle ($\delta_{TG} = 0.715$) }  \\ \cline{2-6}
			$k=2$     & 12 & 13 & 14& 15 & 15 \\
			$k=3$     & 11 & 11 & 11 & 11 & 11 \\
			$k=4$     & 11 & 11 & 11 & 11 & 11 \\
			\hline
			&\multicolumn{5}{c}{C-AMLI-cycle ($\delta_{TG} =  1$)}  \\ \cline{2-6}
			$k=2$     & 12 & 13 & 13 & 14 & 15 \\
			$k=3$     & 14 & 14 & 14 & 14 & 14 \\
			$k=4$     & 24 & 25 & 30 & 40 & 45 \\
			\hline
			&\multicolumn{5}{c}{M-AMLI-cycle }  \\ \cline{2-6}
			$k=2$     & 12 & 13 & 13 & 14 & 15 \\
			$k=3$     & 11 & 11 & 11 & 11 & 11 \\
			$k=4$     & 11 & 11 & 11 & 11 & 11 \\
			\hline  \hline
		\end{tabular}
	\end{center}
\end{table}

\begin{table}[H]
	\begin{center}
		\caption{CPU time for Example~\ref{exp:aniso-poisson}.} 
		\label{tab:Exa.aniso-poissontwogrid and cycles-CPUtime}
		\renewcommand{\arraystretch}{1.0}
		\begin{tabular}{cccccc}
			\hline \hline
			&  $h=1/128$ & $h=1/256$ & $h=1/512$ & $h=1/1024$ & $h=1/2048$ \\ \hline
			&\multicolumn{5}{c}{K-cycle }  \\ \cline{2-6}
			$k=2$     & 0.072 & 0.190 & 0.806  & 2.693   & 10.482 \\
			$k=3$     & 0.097 & 0.309 & 1.179  & 4.131  & 16.112 \\
			$k=4$     & 0.143 & 0.570 & 2.075 & 7.067  & 28.973 \\
			\hline
			&\multicolumn{5}{c}{C-AMLI-cycle ($\delta_{TG} = 0.715$) }  \\ \cline{2-6}
			$k=2$     & 0.029 & 0.126 & 0.462 & 1.890 & 7.221 \\
			$k=3$     & 0.048 & 0.198 & 0.671 & 2.416 & 9.347 \\
			$k=4$     & 0.112  & 0.386 & 1.243 & 4.298 & 17.570 \\
			\hline
			&\multicolumn{5}{c}{M-AMLI-cycle }  \\ \cline{2-6}
			$k=2$     & 0.027 & 0.149 & 0.423  & 1.772   & 7.157 \\
			$k=3$     & 0.051 & 0.197 & 0.658  & 2.359  & 9.302 \\
			$k=4$     & 0.088 & 0.352 & 1.264 & 4.434 & 16.877  \\
			\hline  \hline
		\end{tabular}
	\end{center}
\end{table}

In~Table~\ref{tab:Exa.aniso-poissontwogrid and cycles}, we present the number of iterations of PCG with different MG cycle preconditioners for different mesh sizes. The overall results are consistent with those observed in Example~\ref{exp:poisson}. First, we can see that the two-grid method converges uniformly. However, the performance of the V-cycle (i.e., $k$V-cycle with~$k=1$) and the W-cycle (i.e., $k$V-cycle with~$k=2$) degenerates as before. The~$k$V-cycle achieves uniform convergence when~$k = 4$. The K-cycle achieves uniform convergence for all $k$, confirming its expected optimality. Both the C-AMLI-cycle and M-AMLI-cycle perform similarly in this case. Specifically, the numbers of iterations for both cycles increase slightly for $k=2$ and remain uniform for $k=3$, $4$, with $11$ iterations across all cases, which is comparable with the K-cycle. On the other hand, if we do not estimate $\delta_{TG}$ (i.e., use $\delta_{TG} = 1$) for C-AMLI-cycle,  we do see the number of iterations grows as $k$ increases for a fixed $h$. Furthermore, the number of iterations increases as $h$ gets smaller for $k=4$. Again, these observations confirm the theoretical results presented in  Section~\ref{sec:AMLI-Cycle-Chebyshev} and Section~\ref{sec:AMLI-Cycles-MA}.

In addition, we compare the CPU times of the K-cycle, C-AMLI-cycle (with $\delta_{TG} = 0.715$), and M-AMLI-cycle in Table~\ref{tab:Exa.aniso-poissontwogrid and cycles-CPUtime}. Similar to the previous results, we observe that both the C-AMLI-cycle and M-AMLI-cycle are approximately twice as fast as the K-cycle. This further confirms our expectations and suggests that the M-AMLI-cycle is a preferable choice in practice.

\subsection{Jump Coefficient Problem}
\begin{figure}[h!]
	\centering
	\includegraphics[scale=0.8]{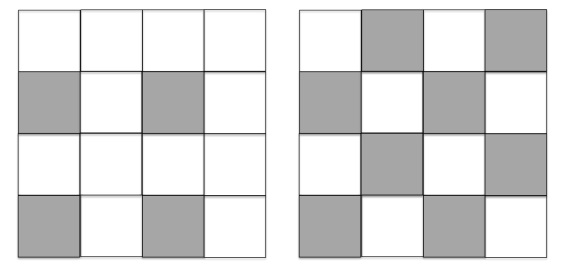}
	\caption{Distribution of the jump coefficient $a(\mathbf{x})$. Left: P1. Right: P2. }
	\label{fig:example2}
\end{figure}

The model problem in this section is an elliptic diffusion problem with jump coefficients. 
\begin{example}\label{exp:jump-poisson}
	Let $\Omega = [0,1] \times [0,1]$, consider
	\begin{equation*}
		\begin{array}{rcl}
			- \nabla \cdot (a(\mathbf{x}) \nabla \mathbf{u}) &=& f, \ \ \ \text{in} \ \ \Omega ,\\
			\mathbf{u} &=& 0, \ \ \ \text{on} \ \ \partial \Omega.
		\end{array}
	\end{equation*}
\end{example}
We consider the case where the diffusion coefficients $a(\mathbf{x})$ are highly oscillatory. Two different test problems are considered, corresponding to different distributions of the diffusion jump $a(\mathbf{x})$. We comment that those two cases were considered in \cite{brannick2018optimal}. In the first problem P1,  the interfaces of the jumps do not intersect, namely, 
\begin{equation*}
	a(\mathbf{x})=
	\begin{cases}
		1, &  \mathbf{x} \in \Omega_1, \\
		10^{-k_{ij}}, &  \mathbf{x} \in  \Omega \backslash \Omega_1,
	\end{cases}
\end{equation*}
where the domain~$ \Omega_1$~corresponds to the white regions in the plot on the left in ~Figure~\ref{fig:example2}. The values of~$k_{ij}\in\{1,2,3,...,6\}$~are selected randomly with a uniform distribution (using built-in MATLAB function \texttt{randi}). In the second problem P2, we consider a checkerboard pattern for the distribution of the jumps, where $\Omega_1$ now corresponds to the white regions in the plot on the right in ~Figure~\ref{fig:example2}. For P2, we randomly select the values~$k_{ij}$ as in P1.

We use a standard cell-centered finite-volume method (see \cite{eymard2000finite,samarskii2001theory}) to discretize P1 and P2 on a structure grid of $\Omega$.  We refer to \cite{brannick2018optimal} for more details.  As before, we set $f = 0$ for both P1 and P2, which results in a zero exact solution.   Similar to the anisotropic case Example~\ref{exp:aniso-poisson}, we set the strength connection parameter to be $0.25$ for dropping small entries in the matrices and forming aggregation accordingly to construct the hierarchy of the UA-AMG method.   When $h = 1/2048$,  the minimal coarsening ratio is about $2$ for both P1 and P2, while the averages are $4.3226$ and $3.7023$, respectively.  Therefore,  we use $k\leq 3$ for Example~\ref{exp:jump-poisson} to maintain nearly optimal computational complexity.

\begin{table}[H]
	\begin{center}
		\caption{Number of PCG iterations for Example~\ref{exp:jump-poisson} P1.}
		\label{tab:Exa.jump-poissontwogrid and cycles-P1}
		\renewcommand{\arraystretch}{1.0}
		\begin{tabular}{ccccccc}
			\hline \hline
			&  $h=1/128$ & $h=1/256$ & $h=1/512$ & $h=1/1024$ & $h=1/2048$\\ \hline
			&\multicolumn{5}{c}{Two-grid}  \\  \cline{2-6}
			& 11 & 11 & 11 & 11 & 11 \\ 
			\hline
			&\multicolumn{5}{c}{ $k$V-cycle  }  \\ \cline{2-6}
			$k=1$   & 23 & 33 & 50 & 66& 86\\
			$k=2$   & 14 & 15 & 16 & 17 & 17 \\
			$k=3$   & 12 & 12 & 12 & 12 & 12 \\
			\hline
			&\multicolumn{5}{c}{K-cycle }  \\ \cline{2-6}
			$k=2$   & 11 & 12 & 12 & 12 & 12 \\
			$k=3$   & 11 & 11 & 11 & 11 & 11 \\
			\hline
			&\multicolumn{5}{c}{C-AMLI-cycle  ($\delta_{TG} = 0.745$) }  \\ \cline{2-6}
			$k=2$   & 12 & 12 & 12 & 12 & 12 \\
			$k=3$   & 11 & 12 & 12 & 12 & 12 \\
			\hline
			&\multicolumn{5}{c}{C-AMLI-cycle ($\delta_{TG} = 1$) }  \\ \cline{2-6}
			$k=2$   & 12 & 12  & 12  & 12  &  12   \\
			$k=3$   & 15 & 15 & 15 &   15  &  15   \\
			\hline
			&\multicolumn{5}{c}{M-AMLI-cycle }  \\ \cline{2-6}
			$k=2$    & 12 & 12 & 12 & 12 & 12 \\
			$k=3$    & 11 & 11 & 11 & 11 & 11 \\
			\hline  \hline
		\end{tabular}
	\end{center}
\end{table}

\begin{table}[H]
	\begin{center}
		\caption{CPU time for Example~\ref{exp:jump-poisson} P1.}
		\label{tab:Exa.jump-poissontwogrid and cycles-P1-CPUtime}
		\renewcommand{\arraystretch}{1.0}
		\begin{tabular}{ccccccc}
			\hline \hline
			&  $h=1/128$ & $h=1/256$ & $h=1/512$ & $h=1/1024$ & $h=1/2048$\\ \hline
			&\multicolumn{5}{c}{K-cycle }  \\ \cline{2-6}
			$k=2$   & 0.097 & 0.446  & 1.605 & 5.981 & 22.025 \\
			$k=3$   & 0.225  & 1.014  & 4.101  & 15.542  & 57.586 \\ 
			\hline
			&\multicolumn{5}{c}{C-AMLI-cycle  ($\delta_{TG} = 0.745$) }  \\ \cline{2-6}
			$k=2$   & 0.060 & 0.244 & 0.953 & 3.487 & 13.124 \\
			$k=3$   & 0.134 & 0.645 & 2.738 & 10.025 & 35.962 \\
			\hline
			&\multicolumn{5}{c}{M-AMLI-cycle }  \\ \cline{2-6}
			$k=2$    & 0.061  & 0.251  & 0.921  & 3.381  & 12.610 \\
			$k=3$    & 0.141 & 0.576 & 2.814 & 8.969 & 32.314 \\
			\hline  \hline
		\end{tabular}
	\end{center}
\end{table}

\begin{table}[H]
	\begin{center}
		\caption{Number of PCG iterations for Example~\ref{exp:jump-poisson} P2}
		\label{tab:Exa.jump-poissontwogrid and cycles-P2}
		\renewcommand{\arraystretch}{1.0}
		\begin{tabular}{cccccc}
			\hline \hline
			&  $h=1/128$ & $h=1/256$ & $h=1/512$& $h=1/1024$ &  $h=1/2048$  \\ \hline
			&\multicolumn{5}{c}{Two-grid}  \\  \cline{2-6}
			& 12 & 11 & 12 & 12 & 12    \\ 
			\hline
			&\multicolumn{5}{c}{ $k$V-cycle  }  \\ \cline{2-6}
			$k=1$    & 64 & 76 & 106 & 146& 194\\
			$k=2$    & 31 & 29 & 37  &  34 & 36\\
			$k=3$    & 23 & 23 & 29 &  26 & 26 \\
			\hline
			&\multicolumn{5}{c}{K-cycle }  \\ \cline{2-6}
			$k=2$    & 20 & 21 & 23 & 26 & 24 \\
			$k=3$    & 14 & 16 & 17 & 19 & 18\\
			\hline
			&\multicolumn{5}{c}{C-AMLI-cycle ($\delta_{TG}=0.749$) }  \\ \cline{2-6}
			$k=2$    & 23 & 27 & 26 & 30 & 30 \\
			$k=3$    & 16 & 17  & 19 & 21  & 20 \\
			\hline
			&\multicolumn{5}{c}{C-AMLI-cycle ($\delta_{TG}=1$) }  \\ \cline{2-6}
			$k=2$    & 22 & 30 & 31 & 31 & 32 \\
			$k=3$    & 22 & 27 & 45 & 40 & 39\\
			\hline
			&\multicolumn{5}{c}{M-AMLI-cycle }  \\ \cline{2-6}
			$k=2$    & 23 & 27 & 26 & 29 & 30\\
			$k=3$    & 17 & 18 & 20 & 22 & 21 \\
			\hline  \hline
		\end{tabular}
	\end{center}
\end{table}

\begin{table}[H]
	\begin{center}
		\caption{CPU time for Example~\ref{exp:jump-poisson} P2}
		\label{tab:Exa.jump-poissontwogrid and cycles-P2-CPUtime}
		\renewcommand{\arraystretch}{1.0}
		\begin{tabular}{cccccc}
			\hline \hline
			&  $h=1/128$ & $h=1/256$ & $h=1/512$& $h=1/1024$ &  $h=1/2048$  \\ \hline
			&\multicolumn{5}{c}{K-cycle }  \\ \cline{2-6}
			$k=2$    & 0.174  & 0.800  & 3.296 &  14.399  & 52.027 \\
			$k=3$    & 0.280 & 1.411    & 5.817  &   28.512  & 100.392 \\
			\hline
			&\multicolumn{5}{c}{C-AMLI-cycle ($\delta_{TG}=0.749$) }  \\ \cline{2-6}
			$k=2$    & 0.119 & 0.576 & 2.009  & 8.883  & 35.090 \\
			$k=3$    & 0.216 & 0.960  & 4.114 & 19.542 & 67.573   \\
			\hline
			&\multicolumn{5}{c}{M-AMLI-cycle }  \\ \cline{2-6}
			$k=2$    & 0.123 & 0.553  & 1.956  & 8.265 & 33.116\\
			$k=3$    & 0.222 & 0.943 & 4.137& 19.899 & 69.548 \\
			\hline  \hline
		\end{tabular}
	\end{center}
\end{table}

We present the number of iterations of PCG for Example~\ref{exp:jump-poisson} P1 in~Table~\ref{tab:Exa.jump-poissontwogrid and cycles-P1}.  The two-grid method converges uniformly, as does the K-cycle for all values of $k$. The performance of the V-cycle (i.e., the $k$V-cycle with $k=1$) degenerates as expected. When $k=3$, the $k$V-cycle becomes uniform.   In this example, for the mesh with $h=1/128$, the computed two-grid convergence rate is $0.7421$, and thus, we use $\delta_{TG} = 0.745$ for P1. Both proposed C-AMLI-cycle and M-AMLI-cycle exhibit uniform convergence for all $k$, demonstrating their effectiveness. However, without estimating $\delta_{TG}$ (just set it to be $1$), the number of iterations increases with $k$, although it is still uniform with respect to $h$. 

As expected, Example~\ref{exp:jump-poisson} P2 is the most challenging case. On the mesh with $h=1/128$, the convergence rate of the two-grid method is $0.7685$, which exceeds $0.75$ and is insufficient for C-AMLI-cycle to converge uniformly for $k=2$. Therefore, to ensure the C-AMLI-cycle is well-defined, we use $\delta_{TG} = 0.749$ in our implementation for P2.   

The number of iterations of PCG for Example~\ref{exp:jump-poisson} P2 are reported in ~Table~\ref{tab:Exa.jump-poissontwogrid and cycles-P2}. As expected, the two-grid method converges uniformly, while the performance of the V-cycle (i.e., $k$V-cycle with $k=1$) degenerates. In fact, the~$k$V-cycle still does not converge uniformly with respect to~$h$ when~$k=2$ or $k=3$. In contrast, the K-cycle is nearly uniform, with the number of iterations increases slightly for $k=2$ and $k=3$. 
The general behavior of the C-AMLI-cycle and M-AMLI-cycle methods is similar, with a slight increase in the number of iterations for both C-AMLI-cycle with estimating $\delta_{TG}$ and M-AMLI-cycle. Nevertheless, the growth is comparable with the K-cycle and both cycles require just a few more iterations compared with the K-cycle. Overall, their performance is still nearly optimal.  On the other hand, when $\delta_{TG}$ is not estimated for the C-AMLI-cycle, it fails to converge uniformly with respect to $h$, and its performance worsens as $k$ increases.

Finally, Tables~\ref{tab:Exa.jump-poissontwogrid and cycles-P1-CPUtime} and~\ref{tab:Exa.jump-poissontwogrid and cycles-P2-CPUtime} present the CPU times of the K-cycle, C-AMLI-cycle, and M-AMLI-cycle for Example~\ref{exp:jump-poisson} P1 and P2, respectively. As before, both the C-AMLI-cycle and M-AMLI-cycle are faster than the K-cycle, as expected. Moreover, the results highlight the practical advantage of the M-AMLI-cycle, even in more challenging applications.

\section{Conclusions}\label{sec:conc}
In this work, we revisit and enhance the AMLI-cycle by simplifying both its theoretical foundation and practical implementation. By leveraging the properties of the Chebyshev polynomials, we remove the need to estimate extreme eigenvalues at all levels. Instead, only the two-grid convergence rate needs to be estimated, and this can be computed on a coarse level in practice. This reduces computational complexity, making the AMLI-cycle more practical for large-scale problems.  Additionally, we introduce a momentum-accelerated AMLI-cycle, drawing inspiration from optimization techniques such as Nesterov acceleration and stationary Anderson Acceleration. This new approach achieves a  uniformly bounded condition number without estimating the extreme eigenvalues or two-grid convergence rate, simplifying its implementation to be as straightforward as standard multigrid methods. Our theoretical analysis shows that the momentum-accelerated AMLI-cycle is asymptotically as good as the Chebyshev-based AMLI-cycle for $k=2$ case and nearly optimal for $k \ge 3$. In our preliminary numerical tests, both the proposed Chebyshev-based and momentum-accelerated AMLI-cycles remain robust and exhibit iteration counts comparable to those of the K-cycle across various scenarios. Moreover, both variants are approximately twice as fast as the K-cycle in terms of CPU time. These results reinforce the potential of the AMLI-cycle as a practical and efficient approach for a wide range of applications.

For future work,  we would like to investigate other polynomials so that the resulting AMLI-cycle could be asymptotically as good as the Chebyshev-based AMLI-cycle for any $k$ without the need of estimating extreme eigenvalues and two-grid convergence rate.   In addition, since momentum accelerations can be applied to general nonconvex optimization problems, we plan to extend the momentum-accelerated AMLI-cycle methods for solving general non-SPD linear systems and investigate the performance theoretically and numerically.

\section*{Acknowledgments} The authors wish to thank Ludmil Zikatanov for many insightful and helpful discussions and suggestions. 

Funding: The work of Hu was partially supported by the National Science Foundation [grant number DMS-220826]. The work of Niu was supported by the Natural Science Foundation of Henan Province [grant number 242300421211].


\appendix
\section{Proof of Lemma~\ref{lem:pk<1}} \label{appendix:proof-pk<1}
For a fixed $\widetilde{x} \in (0,1]$, the three-term recurrence relationship \eqref{def:tilde-pk} generates a sequence $\{\widetilde{p}_0, \widetilde{p}_1, \widetilde{p}_2, \cdots \}$ satisfying $\widetilde{p}_{k} = 2 h \widetilde{p}_{k-1} - h \widetilde{p}_{k-2}$ with $h := 1-a\widetilde{x} < 1$.  We rewrite this in the following matrix form, which is a two-term recurrence relationship,
\begin{equation}\label{def:matrix-two-term}
	\begin{pmatrix}
		\widetilde{p}_{k+1} \\ 
		\widetilde{p}_{k}
	\end{pmatrix}
	=
	\begin{pmatrix}
		2h & -h \\
		1  &  0
	\end{pmatrix}
	\begin{pmatrix}
		\widetilde{p}_{k} \\ 
		\widetilde{p}_{k-1}
	\end{pmatrix}.
\end{equation}
The eigenpairs of the $2 \times 2$ matrix above are
$\lambda_1 = h + \sqrt{h^2 - h}, \ \mathbf{v}_1 = \begin{pmatrix}
	\lambda_1 \\
	1  
\end{pmatrix}
\ \text{and} \ \lambda_2 = h - \sqrt{h^2 - h}, \ \mathbf{v}_2 = \begin{pmatrix}
	\lambda_2 \\
	1   
\end{pmatrix}$.
When $\widetilde{x} = \frac{1}{a}$ for $a > 1$, we have the degenerate case that $h=0$ and, thus, $\lambda_1 = \lambda_2 = 0$. In this case, $\widetilde{p}_k = 0$, $k \geq 2$. Therefore, $\widetilde{x}=\frac{1}{a}$ becomes an root which, of course, satisfies $|\widetilde{p}_k(\widetilde{x})| < 1$.  Otherwise, $\lambda_1 \neq \lambda_2$, thus $\mathbf{v}_1$ and $\mathbf{v}_2$ form a basis of $\mathbb{R}^2$ and we have, 
\begin{equation} \label{eqn:p0-p1-eigen-decomp}
	\begin{pmatrix}
		\widetilde{p}_{1} \\ 
		\widetilde{p}_{0}
	\end{pmatrix}
	=\begin{pmatrix}
		1-\widetilde{x}\\ 
		1
	\end{pmatrix}
	= c \mathbf{v}_1+ (1-c) \mathbf{v}_2 \ \text{with} \ c=\frac{1-\widetilde{x}-\lambda_2}{\lambda_1-\lambda_2}.
\end{equation}
From~\eqref{def:matrix-two-term} and \eqref{eqn:p0-p1-eigen-decomp}, we have
$\begin{pmatrix}
	\widetilde{p}_{k} \\ 
	\widetilde{p}_{k-1}
\end{pmatrix}
= 	\begin{pmatrix}
	2h & -h \\
	1  &  0
\end{pmatrix}^{k-1}
\begin{pmatrix}
	\widetilde{p}_1 \\ 
	\widetilde{p}_0
\end{pmatrix}
=c\lambda_1^{k-1} \mathbf{v}_1+ (1-c)\lambda_2^{k-1} \mathbf{v}_2$,
which implies 
\begin{equation} \label{eqn:p_k+1}
	\widetilde{p}_{k}=c\lambda_1^{k}+ (1-c)\lambda_2^{k}.
\end{equation}
Next we consider different cases based on \eqref{eqn:p_k+1}.

{\bf Case 1:} If $a=1$, then $h=1-\widetilde{x} \in (0,1)$ for $\widetilde{x} \in (0,1)$ and, thus $c=\frac{h-(h-\sqrt{h^2-h})}{2\sqrt{h^2-h}}=\frac{1}{2}$.  In this case, $\lambda_1$ and $\lambda_2$ are complex conjugate and
$
|\lambda_1| = |\lambda_2|  = \sqrt{ h^2 +  (\sqrt{h-h^2})^2 } = \sqrt{h} < 1.
$
Thus, from \eqref{eqn:p_k+1}, we have
$
|\widetilde{p}_{k}|\leq (|\lambda_1|^{k}+|\lambda_2|^{k})/2<1.
$
When $\widetilde{x}=1$, then we have $h=0$, which is the degenerate case and $|\widetilde{p}_k|<1$ holds as discussed before.

{\bf Case 2:} If $\frac{1}{2} \leq a<1$ for $\widetilde{x}\in(0,1]$ or $1<a<\frac{4}{3}$ for $0<\widetilde{x}<\frac{1}{a} < 1$, then $ h = 1-a\widetilde{x} \in (0,1)$.  In this case,  the two roots are complex conjugate and $|\lambda_{1}| = |\lambda_2|=\sqrt{h}<1$.   Moreover, 
$c = \frac{1-\widetilde{x} -\lambda_2}{\lambda_1 - \lambda_2} = \frac{1}{2} - \frac{1-\widetilde{x}-h}{2 \sqrt{h-h^2}} i$ and $  1-c = \overline{c}$.
Thus, from \eqref{eqn:p_k+1}, we have
$
\widetilde{p}_{k} = c \lambda_1^k + \overline{c} \overline{\lambda_1}^k = c \lambda_1^k + \overline{c \lambda_1^k} = 2 \mathfrak{Re}(c \lambda_1^k).
$
Therefore, 
$|\widetilde{p}_k| \leq 2 |c| |\lambda_1|^k = 2 \sqrt{\frac{1}{4} + \frac{(1-\widetilde{x}-h)^2}{4(h-h^2)}} \sqrt{h}^k = \sqrt{h^k + \frac{(1-\widetilde{x}-h)^2}{(1-h)} h^{k-1}}$.
Note that, when $a \geq \frac{1}{2}$, $(1-\widetilde{x}-h)^2 = (a-1)^2\widetilde{x}^2 \leq  a^2 \widetilde{x}^2 = (1-h)^2$ since $(a-1)^2 \leq a^2$ when $a \geq \frac{1}{2}$.  Thus, we obtain that
$
|\widetilde{p}_k| \leq \sqrt{h^k + (1-h) h^{k-1}} = \sqrt{h^{k-1}} < 1.
$

{\bf Case 3:} If $1<a<\frac{4}{3}$ and $\frac{1}{a}\leq \widetilde{x} \leq 1$, then $h=0$ when $\widetilde{x} = \frac{1}{a}$, which is the degenerate case discussed before and $|\widetilde{p}_k|<1$ holds.  For $\frac{1}{a} < \widetilde{x} \leq 1$, $1-a < h < 0$ and, thus, $\lambda_1$ and $\lambda_2$ are real.  Note that, when $1< a < \frac{4}{3}$, we have 
$
h = 1-a\widetilde{x} \geq 1 - a > -\frac{1}{3}. 
$
Thus,  
$
|\lambda_1| < |\lambda_2| = |h - \sqrt{h^2 - h}| < 1.
$
Moreover, $\lambda_1 > 0 > \lambda_2$. It follows that $\lambda_1-\lambda_2>0, 1-\widetilde{x}-\lambda_2>0$, and $
c=\frac{1-\widetilde{x}-\lambda_2}{\lambda_1-\lambda_2}>0$.  If $c\geq1$, then
$|\widetilde{p}_k| \leq |c \lambda_1| + |(1-c)\lambda_2| =c\lambda_1 + (1-c)\lambda_2=1-\widetilde{x}<1$.
On the other hand, if $0<c<1$, then
$	|\widetilde{p}_k| \leq |c \lambda_1| + |(1-c)\lambda_2| \leq c |\lambda_1| + (1-c) |\lambda_2| < c + (1-c) = 1$.

{\bf Case 4:} If $a=\frac{4}{3}$ and $\frac{3}{4} \leq \widetilde{x} \leq 1$, then $h = 0$ when $\widetilde{x} = \frac{3}{4}$ and $h > -\frac{1}{3}$ when $\frac{3}{4} < \widetilde{x} <1$.  Following the same discussion in Case 3, we have $|\widetilde{p}_k|<1$.  Now, when $\widetilde{x} = 1$,  we have $h = -\frac{1}{3}$. This implies $\lambda_1 = \frac{1}{3}$, $\lambda_2 = -1$, and $c = \frac{3}{4}$.   Thus,
$|\widetilde{p}_k| = \left| \frac{3}{4} \left( \frac{1}{3} \right)^k + \frac{1}{4} \left( -1 \right)^k  \right| < 1$.
By summarizing all the cases, we complete the proof.

\end{document}